\documentclass{article}

\usepackage{geometry}               
\usepackage{graphicx}
\usepackage[usenames,dvipsnames]{xcolor}
\usepackage{amssymb,amsmath,amsthm,amsfonts}
\usepackage{bm}                     
\usepackage[english]{babel}
\usepackage[utf8]{inputenc}
\usepackage[colorlinks]{hyperref}
\hypersetup{linkcolor=blue,citecolor=blue,filecolor=black,urlcolor=blue}

\usepackage{enumitem}

\newtheorem{proposition}{Proposition}[section]
\newtheorem{theorem}[proposition]{Theorem}

\newtheorem{lemma}[proposition]{Lemma}

\newtheorem{remark}[proposition]{Remark}
\theoremstyle{definition}

\numberwithin{equation}{section}
\newcommand{\abs}[1]{\left|#1\right|}
\newcommand{\norm}[1]{\left\lVert#1\right\rVert}

\newcommand{\eps}{\varepsilon}

\newcommand{\N}{{\mathbb{N}}}
\newcommand{\Z}{{\mathbb{Z}}}

\newcommand{\R}{{\mathbb{R}}}
\newcommand{\C}{{\mathbb{C}}}

\newcommand{\sphere}{{\mathbb{S}}}

\newcommand{\n}{{n}}

\newcommand{\Ader}[2]{\nabla_{#2}#1}
\newcommand{\AderD}[2]{\mathcal{D}^{3}_{#2}#1}
\newcommand{\tor}{\mathbb{T}^2}

\newcommand{\Bext}{B_{\rm ext}}
\newcommand{\bext}{b_{\rm ext}}

\newcommand{\Hcal}{{\mathcal{H}}}

\DeclareMathOperator{\supp}{supp}

\DeclareMathOperator{\diverg}{div}

\DeclareMathOperator{\real}{Re}
\DeclareMathOperator{\immaginary}{Im}

\usepackage[utf8]{inputenc}
\usepackage[maxbibnames=9, giveninits=true, doi=false, url=false, isbn=false]{biblatex}
\usepackage{doi}
\usepackage{mathrsfs}
\title{A $\Gamma$-convergence result for 2D type-I superconductors }
\author{Alessandro Cosenza, Michael Goldman, Alessandro Zilio }

\begin{document}

\maketitle
\begin{abstract}
    We consider a 2D non-standard Modica-Mortola type functional. This functional arises from the Ginzburg-Landau theory of type-I superconductors in the case of an infinitely long sample and in the regime of comparable penetration  and coherence lengthes. We prove that the functional $\Gamma$-converges to the perimeter functional. This result is a first step in understanding how to extend the results of \cite{CGOS} to the regime of non vanishing Ginzburg-Landau parameter $\kappa$.
\end{abstract}
\noindent
{\footnotesize \textbf{AMS-Subject Classification}}. 
{\footnotesize 35Q56, 49J10, 49J45, 82D55}

\section{Introduction}
Superconductivity was discovered in 1911 by K. Onnes, who observed that certain materials at low temperature showed complete loss of resistivity. In 1933 W. Meissner discovered that this loss of resistivity was related to the expulsion of the magnetic field outside the sample. In 1950, V. Ginzburg and L. Landau proposed a phenomenological model for superconductivity. In their model, a superconductor is represented by a sample $\Omega$ and its state is given through an order parameter $u:\Omega\rightarrow \C$, so that $\rho=\abs{u}$ is the density of superconducting electrons. The magnetic field then lives in the normal region, that is the set in which $\rho$ vanishes, as opposed to the superconducting region, $\rho \simeq 1$, in which no magnetic field is present. This phenomenon is the so-called Meissner effect. Moreover, a superconductor is characterized by two physical parameters, $\lambda$ and $\xi$. The parameter $\lambda$ is called penetration length and represents the typical length of penetration of the magnetic field in the superconducting region, whereas the parameter $\xi$ is called coherence length and represents the typical scale at which $\rho$ varies. The number $\kappa=\lambda/\xi$ is called the Ginzburg-Landau parameter, and its value dictates the behavior of the superconductor. If $\kappa<1/\sqrt{2}$ a positive surface energy penalizes the interface between the normal region and the superconducting region.    
These materials are referred to as type-I superconductors. If instead $\kappa>1/\sqrt{2}$,  the surface tension is negative and the magnetic field penetrates the sample through vortices (see \cite{bbh,serfaty,sandierserfaty,panxiang}). These are instead called type-II superconductors. We are interested in the former type. Type-I superconductors, when subject to an external magnetic field, experimentally exhibit the appearance of complex patterns at the surface \cite{pro1,pro3,pro2}. This behavior hints at the presence of branching patterns inside the sample. A first mathematical justification of this fact was given in \cite{ChCoKoOt,ChKoOt}, where scaling laws for a simplified model were proved. These scaling laws were later extended to the full Ginzburg-Landau model in \cite{COS}. In particular, \cite{ChCoKoOt,ChKoOt, COS} distinguish two regimes of parameters, corresponding to uniform and non uniform branching. In the first regime, the magnetic field distributes uniformly at the boundary on the sample, whereas, surprisingly, in the second regime it is energetically favorable for the magnetic field to deviate from uniformity at the boundary. Going beyond the scaling law, in \cite{CGOS} a reduced model was derived (in the sense of $\Gamma$-convergence \cite{braides,dalmaso}). The resulting model is a non standard branched transport problem, where the magnetic field concentrates on a one dimensional tree structure. Minima of this reduced functional were studied in \cite{minimizers2d,dephilippis2023energy,branched}. The result of \cite{CGOS} was proved in the uniform regime and for vanishing applied magnetic field and  $\kappa$.
\\
The goal of this paper is to better understand the limit behavior of type-I superconductors in the case when $\kappa$ and the external magnetic field do not vanish. Recalling the definition of $\kappa$, this is a regime in which penetration length and coherence length are comparable. 
In the full Ginzburg-Landau Model, various lengthscales are involved, see \cite{CGOS}. We consider a simplified 2D setting, which can be understood informally as a blow-up at the scale of coherence and penetration length. In this setting the Ginzburg-Landau model resembles a non-standard Modica-Mortola type functional \cite{modicamortola,modica}. We prove that this functional $\Gamma$-converges to the perimeter functional. Since $\kappa$ does not vanish and due to the vectorial nature of the problem, several difficulties arise in the proof. In particular, the main difficulty is given by the fact that the functional pairs a local coupling between $\rho$ and the magnetic field $B$ with a global coupling between $\theta$, the phase of $u$, and the magnetic potential $A$. 
\\
We stress again that this result diverges from \cite{CGOS,COS}, since the parameter $\kappa$ is not assumed to be small. One may interpret this result as a first step in understanding how to extend the results of \cite{CGOS} to the regime of non-vanishing $\kappa$.
\\
To present the main result of the paper, let us introduce the functional we study (see Section \ref{sec:model}). Consider the torus $\tor=\R^2\backslash \Z^2$. Let $\eps>0$, $\kappa\in (0,1/\sqrt{2})$ and $\alpha=\kappa^{-1}\eps^{-2}$. The model we consider is
\begin{equation}
\label{e:GL2D}
    \mathcal{E}_\eps(u,A)=\int_{\tor}\eps((1-\kappa\sqrt{2})\abs{\Ader{u}{\alpha A}}^{2}+\kappa\sqrt{2}\abs{\AderD{u}{\alpha A}}^2)+\frac{1}{\eps}\left(B-\frac{1}{\sqrt{2}}(1-\abs{u}^{2})\right)^{2}\,dx.
\end{equation}
with the constraint
\begin{equation*}    
\int_{\tor}B\,dx=\frac{b_{\text{ext}}}{\kappa}.
\end{equation*}
Here, $u$ is the complex valued order parameter representing the state of the material, $\rho^2=\abs{u}^2$ is the density of superconducting electrons, $A$ is the vector valued magnetic potential, $B=\nabla \times A$ is the magnetic field in the sample and $\bext$ is the external magnetic field. The quantity 
$\nabla_{A}u=\nabla u-iAu$   
is the covariant derivative of $u$ and 
$    \AderD{u}{A}=(\nabla_A u)_2-i(\nabla_A u)_1$
is the so-called Bogomol’nyi operator. We consider the setting of periodic boundary conditions on the physical values $\rho$ and $B$(see Section \ref{sec:model} for the precise assumptions on the functional spaces we use). The functional \eqref{e:GL2D} can be derived from the full Ginzburg-Landau model in the case of an infinitely long sample  (see Appendix \ref{sec:derivation}) and the parameter $\eps$ represents the typical width of the interfaces between normal and superconducting regions. Notice that when $\rho>0$, if $u=\rho e^{i\theta}$, then
\begin{equation}
\label{e:splittinglocnonloc}
    \abs{\nabla_{\alpha A} u}^2=\abs{\nabla\rho}^2+\rho^2\abs{\nabla \theta-\alpha A}^2.
\end{equation} 
For $\eps\simeq 0$ and $\rho\simeq 1$, the second term in \eqref{e:splittinglocnonloc} forces $A\simeq \nabla(\theta/\alpha)$, so that $B=0$. This means that for $\eps\simeq 0$, the Meissner condition $\rho^2B\simeq 0$ is enforced. Such condition implies that the last term in $\eqref{e:GL2D}$ can be interpreted as a double well potential in $\rho$ and $B$ with the two phases $(\rho,B)=(1,0)$ and $(\rho,B)=(0,1/\sqrt{2})$. Complementarily, the first term in \eqref{e:splittinglocnonloc} penalizes oscillations of $\rho$. Thus, we can interpret the functional as a Modica-Mortola type functional \cite{modicamortola,modica} in $\rho$ and $B$, which converges to the classical perimeter functional \begin{equation*}
    P(\{\rho=0\},\tor),
\end{equation*}
(see for instance \cite{AFP}). In particular, our functional is of vectorial nature,  similarly to generalizations of the standard Modica-Mortola which have been studied in recent years \cite{dalmasoleoni,contifonsecaleoni,barbarapesic,fortuna2023read,conti2006rigidity,davoli2023two,davoli2020two}.
Our main result is the following. 
\begin{theorem}
\label{th:main}
Let $\kappa\in(0,1/\sqrt{2})$, $\bext<\kappa/\sqrt{2}$ and consider $\eps_{n}\rightarrow 0$. Denoting by  $\sigma_\kappa$ the constant \eqref{e:sigmak} defined in Section \ref{sec:liminf} and by $\mathcal{A}_{\eps_n}$ the space of admissible configurations defined in Section \ref{sec:model}, then the following holds:
\begin{enumerate}
    
    \item \label{i:compactness} Let $(u_{n},A_{n})\in \mathcal{A}_{\eps_n}$ be a sequence such that $\sup_n\mathcal{E}_{\eps_n}(u_{n},A_{n})<\infty$.
    Then, up to a subsequence, $\abs{u_{n}}=\rho_{n}\rightarrow \rho$ in $L^{p}(\tor)$ for all $1\leq p<\infty$, with $\rho\in BV(\tor,\{0,1\})$ and $B_n=\nabla\times A_{n}\rightarrow(1-\rho^2)/\sqrt{2}$ in $L^{2}(\tor)$.
    \item \label{i:liminf} Let $(u_{n},A_{n})\in \mathcal{A}_{\eps_n}$ be a sequence such that
    $\abs{u_{n}}=\rho_{n}\rightarrow \rho$ in $L^{p}(\tor)$ for all $1\leq p<\infty$, with $\rho\in BV(\tor,\{0,1\})$ and $B_n =\nabla\times A_{n}\rightarrow (1-\rho^2)/\sqrt{2}$ in $L^{2}(\tor)$.  
    Then 
    \begin{equation*}
        \liminf_{n} \mathcal{E}_{\eps_n}(u_{n},A_{n})\geq \sigma_{\kappa} P(\{\rho=0\},\tor). 
    \end{equation*}
    \item \label{i:limsup} Suppose $\eps_n^{-2}\kappa^{-2}b_{\text{ext}}\in 2\pi\Z$. Let $\rho \in BV(\tor,\{0,1\})$ with $\int_{\tor} (1-\rho^2)/\sqrt{2}\,dx=b_{\text{ext}}/\kappa$. Then there exists $(u_{n},A_{n})\in \mathcal{A}_{\eps_n}$  sequences such that
    $\abs{u_{n}}=\rho_{n}\rightarrow \rho$ in $L^{p}(\tor)$ for all $1\leq p<\infty$, $B_n=\nabla\times A_{n}\rightarrow (1-\rho^2)/\sqrt{2}$ in $L^{2}(\tor)$, $\int_{\tor}B_n\,dx=b_{\text{ext}}/\kappa$ and 
    \begin{equation*}
        \limsup_n \mathcal{E}_{\eps_n}(u_n,A_n)\leq \sigma_\kappa P(\{\rho=0\},\tor).
    \end{equation*}
\end{enumerate}
\end{theorem}
\begin{remark}
    For Theorem \ref{th:main} to be meaningful we assume
\begin{equation*}
    \frac{b_{\text{ext}}}{\kappa}=\int_{\tor}B\,dx<\frac{1}{\sqrt{2}}.
\end{equation*}
Indeed, if $b_{\text{ext}} /\kappa=1/\sqrt{2}$ the trivial minimizer given by $(\rho,B)=(0,b_{\text{ext}}/\kappa)$ is admissible and no phase transition occurs, whereas for $b_{\text{ext}}/\kappa> 1/\sqrt{2}$ the energy is unbounded in $\eps$ for any $B$. 
\end{remark}
\begin{remark}
    The quantization condition $\eps_n^{-2}\kappa^{-2}b_{\text{ext}}\in 2\pi\Z$ is necessary in the construction of the recovery sequence to ensure that the phase $\theta$ of $u$ is well defined. This is a consequence of the fact that the phase in the full Ginzburg-Landau model is naturally quantized (see \cite{CGOS,COS}).
\end{remark}
\begin{remark}
    If we disregard the quantization condition, we may interpret Theorem \ref{th:main} as a classical $\Gamma$-convergence result. We set \begin{equation*}
        X=\left\{(\rho,B)\in L^{\infty}(\tor)\times L^2(\tor), \rho\leq 1, \int_{\tor}B=\bext/\kappa \right\}
    \end{equation*}  
    and define the functionals $F_\eps,F:X\rightarrow [0,+\infty]$ as
    \begin{equation*}
        F_\eps(\rho,B)=\inf\left\{E_{\eps}(u,A), (u,A)\in\mathcal{A}_\eps, \abs{u}=\rho, \nabla\times A=B\right\}
    \end{equation*}
    and 
    \begin{equation*}
        F(\rho,B)=\begin{cases}
            \sigma_\kappa P(\{\rho=0\},\tor) \qquad &\text{ if }\rho\in BV(\tor,\{0,1\}), \ B=1/\sqrt{2}(1-\rho^2),\\
            +\infty \qquad &\text{ otherwise. }
        \end{cases}
    \end{equation*}
    Theorem \ref{th:main} states that $F_{\eps}\overset{\Gamma}{\rightarrow}F$ with respect on the topology on $X$ given by $L^{p}(\tor)$ convergence in $\rho$ for all $p>1$ and $L^2(\tor)$ convegence in $B$.
\end{remark} 
Let us now comment on the strategy of the proof. Regarding compactness (which we essentially take from \cite{CGOS}), we exploit a weak version of the Meissner effect (Proposition \ref{Prop:weakmeissner}) to prove that the functional \eqref{e:GL2D} can be estimated from below by a classical Modica-Mortola type functional in $\rho$. Then, the compactness follows by the standard Modica-Mortola argument for compactness \cite{modica,albertimodicamortola}.

For the $\Gamma$-$\liminf$, we use a blow up procedure \cite{fonseca_muller,blowup}. We define $\sigma_\kappa$ implicitly as a minimization problem on a unit cell where the interface is a vertical segment.

Next, we prove a characterization for $\sigma_\kappa$ which will be useful in the construction of the recovery sequence. This is the main part of the paper and is divided in two steps.

In the first step (Section \ref{sec:dirichlet}), as in typical Modica-Mortola type results, we show that for the minimization problem on the cell $\sigma_\kappa$ can be interpreted as an optimal profile problem $\sigma_{\kappa,\delta}$ with fixed vertical boundary conditions at a distance $\delta$ from the transition (Proposition \ref{prop:sigmaexchange}). The proof of this result is however delicate in our setting. In particular, extra care is necessary because of the global coupling term in \eqref{e:splittinglocnonloc}. Furthermore, we account for the gauge invariance of the functional (see Section \ref{sec:model}) by fixing an appropriate gauge for which the result holds.

In the second step (Section \ref{subsec:periodic}), we prove that, up to a small cost, we can impose horizontal periodic boundary conditions (Proposition \ref{prop:periodic}). In particular, we impose that $u$ is equal to a real profile $U$ near the top and bottom boundary. We actually conjecture that the optimal transition profile itself is one dimensional, in analogy with some previous results \cite{chapman_91,chapman2000}, see Remark \ref{rem:1D}.
Finally, we prove the $\Gamma$-$\limsup$ inequality. By  density  of polyhedral sets in the class of  finite perimeter sets \cite{maggi}, we can suppose that the set $E=\{\rho=0\}$ is polyhedral. Again in analogy with Modica-Mortola, we ``glue" a rescaled competitor $(u_0,A_0)$ of $\sigma_{\kappa}$ on a $\eps_n$-neighborhood of $\partial E$. Then, we globally construct the physical quantities $\rho_n$ and $B_n$ so that they coincide with the wells of the functional $\mathcal{E}_{\eps_n}$ far from $\partial E$.
The main difficulty of the construction is then to appropriately define the phase $\theta_n$ and the potential $A_n$ so that $u_n=\rho_ne^{i\theta_n}$ and $A_n$ match the local values $u_0,A_0$ on an $\eps_n$ neighborhood of $E$. 
In particular the connected components of $E$ act as ``poles" for the phase $\theta_n$, and force the quantization condition $\eps_n^{-2}\kappa^{-2}b_{\text{ext}}\in 2\pi\Z$. 
%
Let us remark that the condition is less and less stringent as $\eps_n\rightarrow 0$.\\
The rest of the paper is organized as follows. In Section \ref{sec:preliminaries}, we introduce the notation and the precise hypothesis on the model. In Section \ref{sec:compactness}, we prove some a priori estimates and the compactness result. In Section \ref{sec:liminf} we prove the $\Gamma$-$\liminf$ result. Section \ref{sec:char} is devoted to characterizing the constant $\sigma_\kappa$. Finally, in Section \ref{sec:limsup} we prove the $\Gamma$-$\limsup$ inequality.

\section{Preliminaries}
\label{sec:preliminaries}
In this section, we fix the notation and precise the assumptions we make on the energy functional.
\subsection{Notation}
\label{sec:notation}
  We use the notation $a \lesssim b$ to indicate $a\leq C b$ for some universal constant $C > 0$ and  $a \lesssim_p b$ to indicate $a\leq C b$ for some constant $C=C(p) > 0$ . For $L,S>0$ we set $Q_{L}=(-L/2, L/2)^{2}$, $Q_{L,S}=(-L/2,L/2)\times (-S/2,S/2)$ and $H^-_{L,S}=(-L/2,0)\times(-S/2,S/2)$. We call $H^{-}=H^{-}_{1/2,1/2}=(-1/2,0)\times(-1/2,1/2)$ and for any direction $\nu\in \mathbb{S}^1$ we call $H^{-}_{\nu}=R_{\nu}H^{-}$ where $R_{\nu}$ is the rotation matrix such that $R_{\nu}e_1 = \nu$, ($e_i$ being the canonical base of $\R^2$). Moreover, let us denote by $Q^{\nu}_{r}(x_0)$ the cube centered in a point $x_{0}\in Q_1$ of scale $r\in (0,1)$ and of direction $\nu\in \sphere^1$, namely $Q_r^\nu(x_0)=x_{0}+r R_{\nu}Q_{1}$. We also abbreviate $Q_r^\nu=Q_r^\nu(0)$ so that $Q_r=Q_r^{e_1}$. 
  We define the $2D$ curl $\nabla\times v = \partial_1 v_2-\partial_2v_1$. For any vector $v \in \R^2$, we define $v^\perp=(-v_2,v_1)$.  We identify periodic functions on $Q_1$ as functions on the torus $\tor=\R^2/(\Z^2)$ equipped with the distance $\abs{\cdot}=\min_{z\in \Z^2}\abs{\cdot-z}$. We denote $BV(\tor)$ the space of periodic bounded variation functions on $\tor$. For a finite perimeter set $\Omega\subset \tor$(see \cite{AFP} for details), we denote by $P(\Omega,\tor)$ its periodic perimeter and by $\partial^*\Omega$ its reduced boundary. We denote by $\mathcal{H}^1$ the periodic Haussdorff measure on $\tor$. For any measure $\mu$, we denote by $\supp \mu$ its support. For any curve $\Gamma\subset \tor$ whose parametrization is  $\gamma:[0,T]\rightarrow \tor$ and for any field $A:\tor\rightarrow \R^2$, we use the notation $\int_{\Gamma}A\cdot\,dl=\int_{0}^TA(\gamma(t))\cdot \dot{\gamma}(t)\,dt$. We always omit relabeling subsequences. 
\subsection{The 2D model}
\label{sec:model}
We detail the functional setting that will be assumed throughout the paper. For $\eps>0$ and $\kappa\in (0,1/\sqrt{2})$ we let $\alpha=\kappa^{-1}\eps^{-2}$. Let $\Omega \subset \R^2$ be a Lipschitz domain (the material sample). We consider a complex function $u : \Omega \to \C$ which represents the state of the material (the order parameter) and $\rho^2=\abs{u}^2$ the density of superconducting electrons, the vector field $A : \Omega \to \R^2$ describing the magnetic vector potential associated to the magnetic field $B=\nabla \times A$ in the sample an $\bext \in \R$ the external magnetic field. We denote by
\begin{equation*}
 \nabla_{A}u=\nabla u-iAu   
\end{equation*}
the covariant derivative of $u$ and by
\begin{equation*}
    \AderD{u}{A}=(\nabla_A u)_2-i(\nabla_A u)_1
\end{equation*}
the so-called Bogomol’nyi operator. The scalar field
\begin{equation}
\label{e:j}
    j=\immaginary\left(u\overline{\nabla_{\alpha A}u}\right)=\real(-i\overline{u}\nabla_{\alpha A}u),
\end{equation}
describes the superconducting current.

We introduce the set of functions
\begin{equation*}
    \mathcal{A}(\Omega)=\left\{(u,A)\in H^1_{\text{loc}}(\Omega;\mathbb{C})\times L^{2}_{\text{loc}}(\Omega;\mathbb{R}^2),\; \rho \leq 1, \; B \in L_{\text{loc}}^2(\Omega)\right\}
\end{equation*}
and denote $\mathcal{A}=\mathcal{A}(Q_1)$. We impose the condition $\rho\leq 1$ since it is always true for minimizers of the functional \eqref{e:GL2D} (see \cite[Lemma 3.7]{COS} for a proof in the 3D case, the proof in our setting is completely analogous).

We say that two configurations $(u,A)$ and $(\hat{u},\hat{A}) \in \mathcal{A}(\Omega)$ are gauge equivalent if there exists $\varphi\in H^1_{\textup{loc}}(\Omega)$ such that 
\begin{equation}
\label{e:gauge}
\begin{cases}
     u=\hat{u}e^{i\varphi}\\
    A=\hat{A}+\frac{1}{\alpha}\nabla \varphi.
\end{cases}
\end{equation}
Notice that the physical quantities $\rho$, $B$, $j$, and the terms $|\nabla_Au|$ and $|\AderD{u}{A}|$ are invariant under the gauge transformation \eqref{e:gauge}. 

We now specify the space of admissible configurations. We define
\begin{equation*}
    \mathcal{A}_{\eps}=\left\{(u,A)\in\mathcal{A}(\R^2) \left| \begin{array}{c} (u(\cdot+k),A(\cdot+k))\text{ is gauge equivalent to }(u,A) \, \forall k\in \Z^2 \\ \int_{Q_1} B \, dx = \bext/\kappa  \end{array}\right.\right\}.
\end{equation*}
Here the condition $\int_{Q_1} B \, dx = \bext/\kappa$ is a consequence of Maxwell law for the magnetic field $B$ (see Appendix~\ref{sec:derivation} for details). Notice that the set depends on $\eps$ since it depends on the gauge equivalence \eqref{e:gauge}, which itself depends on $\alpha=\kappa^{-1}\eps^{-2}$.
In the set $\mathcal{A}_{\eps}$ the physical quantities $\rho$, $B$ and $j$ are $Q_1$-periodic and we can think of them as defined on the torus $\tor$, whereas in general $u$ and $A$ will not be periodic. We also observe that $(u,A)\in\mathcal{A}_\eps$ implies $\rho\in H^1(\tor)$, $B\in L^2(\tor)$ and $j\in L^2(\tor;\R^2)$. 

For any bounded Lipschitz domain $\Omega \subset \R^2$ we introduce an analogue of functional \eqref{e:GL2D}. For $(u,A) \in \mathcal{A}(\Omega)$, we let
\begin{equation*}
    \mathcal{E}_\eps(u,A, \Omega)=\int_{\Omega}\eps\left((1-\kappa\sqrt{2})\abs{\Ader{u}{\alpha A}}^{2}+\kappa\sqrt{2}\abs{\AderD{u}{\alpha A}}^2\right)+\frac{1}{\eps}\left(B-\frac{1}{\sqrt{2}}(1-\abs{u}^{2})\right)^{2}\,dx.
\end{equation*}
Notice that $\mathcal{E}_\eps$ is invariant under change of gauge. For brevity, we denote
\begin{equation*}
    \abs{\nabla^{\kappa}_{A}u}^2=(1-\kappa\sqrt{2})\abs{\Ader{u}{ A}}^{2}+\kappa\sqrt{2}\abs{\AderD{u}{  A}}^2,
\end{equation*} so that 
\begin{equation*}
    \mathcal{E}_\eps(u,A, \Omega)=\int_{\Omega}\eps \abs{\nabla^{\kappa}_{\alpha A}u}^2+\frac{1}{\eps}\left(B-\frac{1}{\sqrt{2}}(1-\abs{u}^{2})\right)^{2}\,dx.
\end{equation*} 

We remark that, when $\rho>0$, letting $u=\rho e^{i\theta}$, we can write 
\begin{equation}
\label{e:Aderformula}
    \abs{\nabla_{\alpha A} u}^2=\abs{\nabla\rho}^2+\rho^2\abs{\nabla \theta-\alpha A}^2.
\end{equation}
Moreover, by a direct computation (see \cite[Lemma 2.1]{COS}), we have
\begin{equation}
\label{e:covariantd3relation}
\abs{\Ader{u}{\alpha A}}^{2}=\abs{\AderD{u}{\alpha A}}^{2}+ \alpha\rho^{2}B+\nabla\times j.
\end{equation}

We are interested in studying the limit behavior as $\eps \to 0$ of the functional $\mathcal{E}_\eps(u,A) = \mathcal{E}_\eps(u,A, Q_1)$ on the set of admissible configurations $\mathcal{A}_\eps$ (see Theorem~\ref{th:main}). 
\section{A priori estimates and compactness}
\label{sec:compactness}
The goal of this section is to derive some a priori estimates and prove the following theorem (Theorem \ref{th:main}-\ref{i:compactness}). The proof is essentially contained in the proof of \cite[Proposition 6.1]{CGOS}, but we repeat it for the reader's convenience. 
\begin{theorem}
\label{th:compactess}
Let $\kappa\in(0,1/\sqrt{2})$, $\bext<\kappa/\sqrt{2}$ and consider $\eps_{n}\rightarrow 0$. Let $(u_{n},A_{n})\in \mathcal{A}_{\eps_n}$ be a sequence such that $\sup_n \mathcal{E}_{\eps_n}(u_{n},A_{n})<\infty$. Then, up to a subsequence, $\abs{u_{n}}=\rho_{n}\rightarrow \rho$ in $L^{q}(\tor)$ for all $1\leq q<\infty$, with $\rho\in BV(\tor,\{0,1\})$ and $B_n=\nabla\times A_{n}\rightarrow (1-\rho^2)/\sqrt{2}$ in $L^{2}(\tor)$.
\end{theorem}
We start by proving a weak version of the Meissner effect. Informally, the Meissner effect states that, for $\eps$ small, $\rho B \simeq 0$, that is the magnetic field $B$ is expelled from the region of superconductivity $\{\rho > 0\}$. 
\begin{proposition}
\label{Prop:weakmeissner}
Let $\Omega\subset Q_1$ be a Lipschitz domain. For any $\phi\in H^{1}(\Omega)\cap L^{\infty}(\Omega)$ and $(u,A)\in\mathcal{A}(\Omega)$, we have 
\[
    \abs{\int_{\Omega}\rho^{2}B\phi\, dx}\lesssim_\kappa \eps\left( \mathcal{E}_\eps(u,A,\Omega)\norm{\phi}_{L^\infty(\Omega)} + \eps^{\frac{1}{2}}\mathcal{E}_\eps(u,A,\Omega)^{\frac{1}{2}}\norm{\nabla \phi}_{L^2(\Omega)}+\eps\abs{\int_{\partial \Omega}\phi j\cdot dl}  \right).
    \]

\end{proposition}
\begin{proof}
The proof is essentially the same as \cite[Lemma 6.2]{CGOS}, but without assuming periodic boundary conditions. We start by noticing that 
\begin{equation*}
    \int_{\Omega}(\nabla\times j)\phi\,dx=\int_{\Omega}\nabla\times (j\phi)\,dx-\int_{\Omega} j^\perp\cdot\nabla\phi\,dx=\int_{\partial \Omega}\phi j\cdot dl-\int_{\Omega} j^\perp\cdot\nabla\phi\,dx.
\end{equation*}
Then, multiplying \eqref{e:covariantd3relation} by $\phi$ and integrating over $\Omega$, we have 
\begin{multline*}
    \abs{\int_{\Omega}\rho^{2}B\phi}\leq \eps^{2}\kappa \int_{\Omega} \left(\abs{\Ader{u}{\alpha A}}^{2}+\abs{\AderD{u}{\alpha A}}^{2}\right)\abs{\phi}\,dx + \eps^{2}\kappa \abs{\int_{\Omega}(\nabla\times j)\phi\,dx}\\
   \lesssim_\kappa \eps \mathcal{E}_\eps(u,A)\norm{\phi}_{L^\infty(\Omega)}+\eps^{2}\int_{\Omega}\abs{j}\abs{\nabla \phi}\,dx +\eps^{2}\abs{\int_{\partial \Omega}\phi j\cdot dl}\\ \lesssim_\kappa \eps\left( \mathcal{E}_\eps(u,A)\norm  {   \phi}_{L^\infty(\Omega)}    + \eps^{\frac{1}{2}}\mathcal{E}_\eps(u,A)^{\frac{1}{2}}\norm{\nabla \phi}_{L^2(\Omega)}+\eps\abs{\int_{\partial \Omega}\phi j\cdot dl} \right),
\end{multline*}
where, in the last step, we have exploited the fact that $\rho\leq 1$ implies $\abs{j}\leq \abs{\nabla_{\alpha A} u}$. 
\end{proof}
We now prove that the energy \eqref{e:GL2D} can be estimated from below by a classical Modica-Mortola type functional.
The following lemma is an adaptation of \cite[Lemma 6.5]{CGOS}, see also \cite{bose_einstein} where a similar idea is used.
\begin{lemma}
\label{lem:mmestimate}
Let $\Omega\subset Q_1$ be a Lipschitz domain, $(u,A)\in\mathcal{A}(\Omega)$ and consider the functions
\begin{equation*}
    W(\rho)=\frac{1}{2}\min\{2\rho^2,1\}(1-\rho^{2})^{2},
\end{equation*}

and \begin{equation}
\label{e:psi}
    \psi=\min\left\{2,\frac{1}{\rho^2}\right\}(1-\rho^{2}).
\end{equation}
Then we have the inequality
\begin{equation}
    \label{e:dwcontrol}
    M_\eps(\rho,\Omega)=\int_{\Omega}\left(\eps\abs{\nabla\rho }^2+\frac{1}{\eps}W (\rho)\right)\,dx\lesssim_\kappa \mathcal{E}_\eps(u,A,\Omega)+\eps\abs{\int_{\partial \Omega}\psi j\cdot dl} .
\end{equation}
\end{lemma}
\begin{proof}
Let us start by noticing that, since 

    \begin{equation*}
        \abs{\nabla_{\alpha A}^\kappa u}^{2}\gtrsim_\kappa \abs{\nabla\rho}^{2},
    \end{equation*}
    it is sufficient to prove 
    \begin{equation*}
        \int_{\Omega}\frac{1}{\eps}W(\rho)\,dx\lesssim_\kappa  \mathcal{E}_\eps(u,A)+\eps\abs{\int_{\partial \Omega}\psi j\cdot dl}.
    \end{equation*}
    Young's inequality gives the relation
    \begin{equation*}
        \frac{1}{2}(1-\rho^{2})^{2}\leq(B-\frac{1}{\sqrt{2}}(1-\rho^{2}))^{2}+\sqrt{2}B(1-\rho^{2}).
    \end{equation*}
    Multiplying both sides by $\min\{2\rho^2,1\}/\eps$, we obtain
    \begin{equation}
    \label{e:modicavsenergy}
        \frac{1}{\eps} W (\rho)\leq \frac{1}{\eps}(B-\frac{1}{\sqrt{2}}(1-\rho^{2}))^{2}+\frac{\sqrt{2}}{\eps}\min\{2\rho^2,1\}B(1-\rho^{2}).
    \end{equation}
    Integrating the previous inequality over $\Omega$, we are left with bounding the last term in \eqref{e:modicavsenergy}. By applying Proposition \ref{Prop:weakmeissner} with $\phi=\psi=\min\{2,1/\rho^2\}(1-\rho^{2})$, we have 
    \begin{equation*}
         \int_{Q_1}\frac{1}{\eps}\rho^2B\psi \lesssim_\kappa \mathcal{E}_\eps(u,A)\norm{\psi}_{L^\infty(\Omega)} + \eps^{\frac{1}{2}}\mathcal{E}_\eps(u,A)^{\frac{1}{2}}\norm{\nabla \psi}_{L^2(\Omega)}+\eps\abs{\int_{\partial \Omega}\psi j\cdot dl} . 
    \end{equation*}
Moreover, direct computations give  $\norm{\psi}_{L^\infty(\Omega)}\lesssim 1$ and  $\norm{\nabla \psi}_{L^2(\Omega)}\lesssim \norm{\nabla \rho}_{L^2(\Omega)}$.  Finally, we have 
\begin{equation*}
    \eps^{\frac{1}{2}}\mathcal{E}_\eps(u,A)^{\frac{1}{2}}\norm{\nabla \psi}_{L^2(\Omega)}\lesssim \mathcal{E}_\eps(u,A)+ \eps \norm{\nabla\rho}_{L^2(\Omega)}^2 \lesssim_\kappa  \mathcal{E}_\eps(u,A),
\end{equation*}
where in the last step
we used again that $\abs{\nabla\rho}^2\lesssim_\kappa  \abs{\nabla_{\alpha A}^\kappa u}^2$.
\end{proof}
We are now ready to prove Theorem \ref{th:compactess}.
\begin{proof}[Proof of Theorem \ref{th:compactess}]
    Let $(u_{n},A_{n})\in \mathcal{A}_{\eps_n}$ be a sequence such that  $\sup_n \mathcal{E}_{\eps_n}(u_{n},A_{n})<\infty $.  We apply Lemma \ref{lem:mmestimate} with $\Omega=Q_1$. Since, for $(u_{n},A_{n})\in \mathcal{A}_{\eps_n}$,  $\psi$ and $j$ are periodic, the boundary term in \eqref{e:dwcontrol} vanishes and we find that
    \begin{equation*}
       \sup_n  M_{\eps_n}(\rho_n,\tor)\lesssim \sup_n \mathcal{E}_{\eps_n}(u_{n},A_{n})<\infty.
    \end{equation*}
     Since $M_\eps$ is a Modica Mortola type functional with a double well potential, by the standard Modica-Mortola argument for compactness (see for instance \cite{albertimodicamortola}) we have, up to a subsequence, $\rho_{n}\rightarrow \rho$ in $L^{p}(\tor)$ for all $1\leq p<\infty$, with $\rho\in BV(\tor,\{0,1\})$. Moreover, $\sup_n \mathcal{E}_{\eps_n}(u_{n},A_{n})<\infty$ implies
    \begin{equation*}
        \int_{\tor} \left(B_n-\frac{1}{\sqrt{2}}(1-\rho_n^{2})\right)^2\,dx \lesssim \eps_n.
    \end{equation*}
    This, paired with the $L^p$ convergence of $\rho_n$ to $\rho$, implies that $B_n\rightarrow (1-\rho^2)/\sqrt{2}$ in $L^2(\tor)$.
    \end{proof}
\section{$\Gamma$-$\liminf$ inequality}
\label{sec:liminf}
Recall that we set $H^{-}=(-1/2,0)\times(-1/2,1/2)$. We introduce the set of sequences
\begin{equation}
\label{e:K}
\mathcal{K}=\left\{(u_n,A_n)_{n\in\N}\subset\mathcal{A} \left| \begin{array}{c} \rho_n \rightarrow \chi_{H^{-}} \text{ in } L^{p}(Q_1) \text{ for } p\geq 1 \\ B_n\rightarrow \frac{1}{\sqrt{2}}(1-\chi_{H^{-}}) \text{ in } L^{2}(Q_1) \end{array} \right.\right\}
\end{equation}
and the constant
\begin{equation}
\label{e:sigmak}
         \sigma_{\kappa}=\inf\left\{\liminf_{n} \mathcal{E}_{\eps_n}(u_n,A_n),\eps_n\rightarrow 0,(u_n,A_n)_{n\in\N}\in\mathcal{K}\right\}. 
 \end{equation}
\begin{remark}
    Let $\eps_n\rightarrow 0$. Consider any smooth function $v:\R\rightarrow [0,1]$ such that $v(-\infty)=1$ and $v(x)=0$ for $x>0$. Define $(u_n,A_n)\in \mathcal{K}$ for $(x_1,x_2)\in Q_1$ as 
    \begin{equation*}
        u_n(x_1,x_2)=v\left(\frac{x_1}{\eps_n}\right)=v_n(x_1), \qquad A_n(x_1,x_2)=\left(0,\frac{(x_1)_+}{\sqrt{2}}\right).
    \end{equation*}
    Notice that $\rho_n B_n=0$ and $\abs{\nabla^\kappa_{\alpha A_n} u_n}^2=\abs{\dot{v}_n}^2$. This means that, using a change of variable,
    \begin{equation*}
    \liminf_{n\to \infty} \mathcal{E}_{\eps_n}(u_n,A_n)=\liminf_{n\to \infty}\int_{-\frac{1}{2}}^0\eps_n\abs{\dot{v}_n}^2+\frac{1}{2\eps_n}(1-v_n^2)^2\,dx_1=\int_{-\infty}^0\abs{\dot{v}}^2+\frac{1}{2}(1-v^2)^2\,dt=G(v).
    \end{equation*}
    Hence 
    \begin{equation*}
    \sigma_\kappa\leq \sigma_0= \inf\{G(v): v(-\infty)=1,  v(0)=0\}.
    \end{equation*}
    The minimization problem $\sigma_0$ is, up to a multiplicative constant, the same problem as the one defined in \cite[Lemma 7.2]{CGOS}, and it admits the explicit solution $v(t)=\tanh(t/\sqrt{2})$. In particular this implies $\sigma_\kappa<\infty$. Moreover, arguing  as in  \cite[Lemma 6.5]{CGOS}, it is possible to prove that $\sigma_\kappa\rightarrow \sigma_0$ as $\kappa\rightarrow 0$.
\end{remark}
\begin{remark}
    Since $\sigma_\kappa<\infty$, using a diagonal argument it is always possible to find a sequence $(u_n,A_n)_{n\in\N}\in\mathcal{K}$ such that $\mathcal{E}_{\eps_n}(u_n,A_n)\rightarrow \sigma_\kappa$. 
\end{remark}
 The goal of this section is to prove the following theorem (Theorem \ref{th:main}-\ref{i:liminf}). 
\begin{theorem}
\label{th:liminf}
Let $\kappa\in(0,\frac{1}{\sqrt{2}})$, $\bext<\kappa/\sqrt{2}$ and consider $\eps_{n}\rightarrow 0$. Let $(u_{n},A_{n})\in \mathcal{A}_{\eps_n}$ be a sequence such that
    $\abs{u_{n}}=\rho_{n}\rightarrow \rho$ in $L^{p}(\tor)$ for all $1\leq p<\infty$, with $\rho\in BV(\tor,\{0,1\})$ and $B_n =\nabla\times A_{n}\rightarrow (1-\rho^2)/\sqrt{2}$ in $L^{2}(\tor)$.
    Then 
    \begin{equation*}
        \liminf_{n} \mathcal{E}_{\eps_n}(u_{n},A_{n})\geq \sigma_{\kappa} P(\{\rho=0\},\tor). 
    \end{equation*}
\end{theorem}
\begin{proof}[Proof of Theorem \ref{th:liminf}]
We use  a blow-up argument.  Let $(u_{n},A_{n})\in \mathcal{A}_{\eps_n}$ be a sequence as in the assumption.
We consider a sequence of positive measures $\{\mu_n\}_{n\in\N}$ on $\tor$ defined on every Borel set $F\subset \tor$ as
\begin{equation*}
    \mu_{n}(F)=\mathcal{E}_{\eps_n}(u_n,A_n,F)=\int_{F} \eps_n\abs{\nabla^{\kappa}_{\alpha A_n}u_n}^2+\frac{1}{\eps_n}\left(B_n-\frac{1}{\sqrt{2}}(1-\abs{u_n}^{2})\right)^{2}\,dx.
\end{equation*}
We may assume without loss of generality that  $\sup_n \mathcal{E}_{\eps_n}(u_n,A_n)<\infty$  and we can extract a subsequence such that $\lim_{n}\mathcal{E}_{\eps_{n}}(u_{\eps_{n}},A_{\eps_{n}})$ exists. By weak compactness of measures (see \cite[Theorem 1.59]{AFP}), up to another subsequence we can suppose that $\mu_{n}$ converges weakly$^{*}$ to a positive finite measure $\mu$.  Now, since $\rho \in BV(\tor,\{0,1\})$, $\rho=\chi_{E}$ for some $E\subset \tor$ with $P(E,\tor) = \Hcal^{1}(\partial^* E)<+\infty$. By the Radon-Nikodym theorem, we can write
\begin{equation*}
    \mu=\frac{d\mu}{d\Hcal^{1}}\Hcal^{1}|_{\partial^* E}+\mu_{c}
\end{equation*}
for some $\mu_{c}\perp \Hcal^{1}|_{\partial^* E}$. We claim that $\frac{d\mu}{d\Hcal^{1}}(x)\geq \sigma_{\kappa}$ for $\Hcal^{1}$-a.e.~$ x\in \partial^* E $. Under such claim, by the semicontinuity of the total variation with respect to weak$^{*}$ convergence (see \cite[Theorem 1.59]{AFP}) we would get
\begin{multline*}
    \liminf_n E_n(u_{\eps_n},A_n)= \liminf_n \mu_{n}(\tor)\geq \mu(\tor)=\int_{\tor}\frac{d\mu}{d\Hcal^{1}}d\Hcal^{1}|_{\partial^* E}+\mu_{c}(\tor)  \\ \geq \int_{\partial^* E}\sigma_{\kappa}d\Hcal^{1}=\sigma_{\kappa}\Hcal^{1}(\partial^* E)=\sigma_{\kappa}P(E,\tor),
\end{multline*}
thus proving the statement. 

Let us now prove the claim. We use the blow-up method (introduced in \cite{fonseca_muller}, see also \cite{blowup}). We know that at $\Hcal^{1}$-a.e.~$x_{0}\in \partial^{*}E$ the set $\partial^{*}E$ admits an outward normal $\nu=\nu(x_{0})$. Using Besicovitch derivation theorem (which we apply using squares instead of balls, see for instance \cite[Theorem 1.153]{fonseca_leoni_calcvar}), the derivative is given by 
\begin{equation*}
    \frac{d\mu}{d\Hcal^{1}|_{\partial^* E}}(x_{0})=\lim_{r\rightarrow 0}\frac{\mu(Q_{r}^{\nu}(x_{0}))}{\Hcal^{1}( Q_{r}^{\nu}(x_{0})\cap \partial^* E )}.
\end{equation*}
We can assume without loss of generality that $x_0=0$. Moreover by the invariance of the energy under rotation we may also assume that $\nu=e_1$. By the properties of sets of finite perimeter we know that (see \cite[Theorems 3.59 and 3.61]{AFP})
\begin{equation*}
    \lim_{r\rightarrow 0}\frac {\Hcal^{1}( Q_r\cap \partial^* E )  }{r}=1,
\end{equation*}
hence 
\begin{equation*}
   \frac{d\mu}{d\Hcal^{1}|_{\partial^* E}}(0)=\lim_{r\rightarrow 0}\frac{\mu(Q_r)}{r}.
\end{equation*}
Since $\mu$ is finite, $\mu(\partial Q_r )=0$ for a.e.~$r$. Thus, by standard weak$^{*}$ convergence arguments (see \cite[Proposition 1.62]{AFP}), we have $\mu_{n}(Q_r)\rightarrow\mu (Q_r)$ for a.e.~$r\in \R$.  Writing for every $\eps_n\rightarrow 0$
\begin{equation*}
    \frac{d\mu}{d\Hcal^{1}|_{\partial^* E}}(0)=\lim_{r\rightarrow 0}\frac{\mu(Q_r)}{\mu_{n} (Q_r)}\frac{\mu_{n} (Q_r)}{r},
\end{equation*}
we can use a diagonal argument to find $\eps_{j},r_{j}$ such that 
\begin{equation}
\label{e:diagderiv}
    \frac{d\mu}{d\Hcal^{1}|_{\partial^* E}}(0)=\lim_{j}\frac{\mu_{j}(Q_{r_{j}})}{r_{j}}.
\end{equation}
Moreover, we can also assume that  $\eps_{j}/r_{j}\rightarrow 0$ and, by the $L^1$ convergence of $\rho_j$,
\begin{equation}
    \label{e:righteps}
    \frac{1}{r_j^2}\int_{\tor}\abs{\rho_j-\chi_E}dx \rightarrow 0.
\end{equation}
By \eqref{e:diagderiv}, we have
\begin{equation*}
    \frac{d\mu}{d\Hcal^{1}|_{\partial^* E}}(0)=\lim_{j}\frac{1}{r_{j}}\int_{Q_{r_{j}}^{\nu}}\eps_{j}\abs{\nabla^{\kappa}_{\alpha_{j}A_{j}}u_{j}}^2+\frac{1}{\eps_{j}}\left(B_{j}-\frac{1}{\sqrt{2}}(1-\abs{u_{j}}^{2})\right)^{2}\,dx.
\end{equation*}
We change variable $x=r_{j}y$ and define $\Hat{A}_{j}(y)=r_{j}^{-1}A_{j}(r_{j}y)$ and $\Hat{u}_{j}(y)=u_{j}(r_{j}y)$ (so that $\hat{\rho}_j=\abs{\hat{u}_j}$ and $\hat{B}_j=\nabla\times \hat{A}_j$). We find
\begin{equation*}
    \frac{d\mu}{d\Hcal^{1}|_{\partial^* E}}(0)=\lim_{j}\int_{Q_{1}}\frac{\eps_{j}}{r_{j}}\abs{\nabla^{\kappa}_{r_{j}^{2}\alpha_{j}\Hat{A}_{j}}\Hat{u}_{j}}^{2}+\frac{r_{j}}{\eps_{j}}\left(\Hat{B}_{j}-\frac{1}{\sqrt{2}}(1-\abs{\Hat{u}_{j}})^{2}\right)^2\,dy.
\end{equation*}
By the properties of sets of finite perimeter (see \cite[Theorem 3.59]{AFP}), we can show that 
for all $p\geq 1, \, \hat{\rho}_j\rightarrow\chi_{H^{-}}$ in $L^{p}(Q_{1})$ and $\Hat{B}_{j}\rightarrow (1-\chi_{H^{-}})/\sqrt{2}$ in $L^{2}(Q_{1})$. Indeed, we have
\begin{multline}
\label{e:convergenceestimate}
\int_{Q_1}\abs{\hat{\rho}_j-\chi_{H^-} }\,dy=\frac{1}{r_j^2}\int_{Q_{r_j}}\abs{\rho_j(x)-\chi_{H^-}\left(\frac{x}{r_j}\right) }\,dx   \\
\leq \frac{1}{r_j^2}\int_{Q_{r_j}}\abs{\rho_{j}(x)-\chi_E(x) }+\abs{\chi_E(x)-\chi_{H^-}\left(\frac{x}{r_j}\right) }\,dx,
\end{multline}    
where the first term on the right of \eqref{e:convergenceestimate} goes to zero by \eqref{e:righteps} and the second one by \cite[Theorem 3.59]{AFP}. Since $\hat{\rho}_j\le 1$ this implies also the convergence in $L^p$ for every $p\ge 1$. Similarly to the proof of Theorem \ref{th:main}-\ref{i:compactness}, the convergence of $\hat{\rho}_j$ and the boundedness of the energy imply the convergence of $\hat{B}_j$.
Hence, setting $\Hat{\eps}_{j}=\eps_{j}/r_{j}$, we find 
\begin{equation*}
    \frac{d\mu}{d\Hcal^{1}|_{\partial^* E}}(0)=\lim_{j}\mathcal{E}_{\Hat{\eps}_{j}}(\Hat{u}_{j},\Hat{A}_{j},Q_1^{\nu}).
\end{equation*}
This concludes the proof of the claimed
\[
 \frac{d\mu}{d\Hcal^{1}|_{\partial^* E}}(0)\ge \sigma_\kappa. \qedhere
\]

\end{proof}
\section{Characterization of $\sigma_\kappa$}
\label{sec:char}
\begin{figure}[h]
    \centering
    \includegraphics[scale = 1]{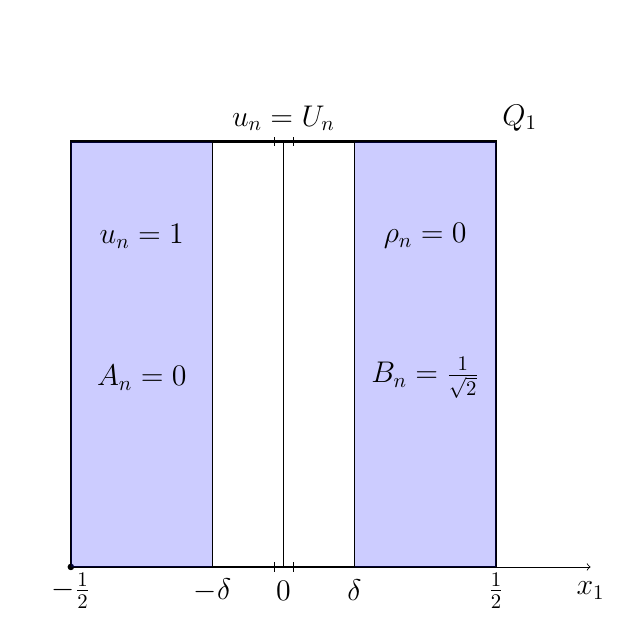} 
    \caption{A typical competitor $(u_n,A_n)_n$ for $\sigma_{\kappa}$. Competitors enjoy fixed boundary conditions at distance $\delta$ from the transition (the zone in blue). Competitors for $\sigma_{\kappa,\delta}^{\text{per
    }}$ moreover are such that $u_n$ coincides with a given real profile $U_n$ near the top and bottom boundaries of $Q_1$.}
    \label{fig:competitor}
\end{figure}

In this section we prove a characterization for $\sigma_\kappa$, which will help with the construction of the recovery sequence. This is done in two steps. First, we show that, at a distance $0<\delta<1/2$ from the interface, we can impose vertical Dirichlet boundary conditions for competitors of $\sigma_\kappa$. In analogy with the study of Modica Mortola type functionals, that means that we think of $\sigma_\kappa$ as an optimal profile problem on the unit cell. Then, we prove that we can also impose horizontal periodic boundary conditions, up to a small cost which depends on $\delta$. That is, we prove that for any $0<\delta<1/2$, we can define a minimization problem $\sigma_{\kappa,\delta}^{\text{per}}$, see \eqref{e:sigmaperdelta}, such that
   \begin{equation*}
         \sigma_{\kappa,\delta}^{\text{per}}\leq \sigma_{\kappa}+C\delta
     \end{equation*}
     and whose competitors enjoy horizontal periodic boundary conditions (see Figure \ref{fig:competitor}).
\subsection{Characterization of $\sigma_\kappa$: vertical Dirichlet conditions}
\label{sec:dirichlet}
Define, for $0<\delta<1/2$, the set of functions
\begin{equation*}
     \mathcal{A}^{\delta}=\left\{(u,A)\in\mathcal{A} \left| \begin{array}{c} (u,A)=(1,0) \text{ for } x_1<-\delta \\ (u,B)=(0,\frac{1}{\sqrt{2}}) \text{ for } x_1>\delta  \\ A_{1}=0 \end{array} \right.\right\},
\end{equation*}
the set of sequences
\begin{equation}
\label{e:BC}
    BC^{\delta}=\left\{(u_\n,A_\n)_\n\in\mathcal{K} \,|\, (u_n,A_n)\in \mathcal{A}^\delta \right\}
\end{equation}
and the constant
    \begin{equation}
    \label{e:sigmabar}
        \sigma_{\kappa,\delta}=\inf\left \{\liminf_{\n}\mathcal{E}_{\eps_n}(u_\n,A_\n) \,|\,\eps_n\rightarrow 0, (u_\n,A_\n)_\n\in BC^\delta\right\}.
     \end{equation}
    The goal of this section is to prove the following proposition. 
\begin{proposition}
\label{prop:sigmaexchange}
For every $0<\delta<1/2$,
\begin{equation*}
\sigma_{\kappa}= \sigma_{\kappa,\delta}.    
\end{equation*}
\end{proposition}
\begin{remark}
    The conditions on $(u_\n,A_\n)$ in \eqref{e:BC} come from two different requirements. On the one hand, as in usual Modica-Mortola type arguments, we want to define the constant $\sigma_{\kappa,\delta}$ as a minimization problem with fixed boundary conditions at some distance $\delta$ from the transition line, as this simplifies the construction of a recovery sequence. This explains the conditions on $\rho_\n$ and $B_\n$. On the other hand, the conditions $u_\n=1$ for $x_1$ near  $-1/2$ and $A_{\n,1}=0$ fix the gauge of the minimizing sequence. This is helpful for adding horizontal periodic boundary conditions (see Section \ref{subsec:periodic}), which in turn helps again with the construction of the recovery sequence.
\end{remark}
Let us start by understanding the scaling properties of the constants $\sigma_{\kappa}$ and $\sigma_{\kappa,\delta}$. First, we make explicit their dependence on the domain. 
Recalling the notation of Section \ref{sec:notation}, for $a,b>0$ we consider the sets
\begin{equation*}
    \mathcal{A}_{a,b}=\{(u,A)\in H^1(Q_{a,b};\mathbb{C})\times L^{2}(Q_{a,b};\mathbb{R}^2) \left| \rho\leq 1, \nabla \times A \in L^2(Q_{a,b}) \right.\}
\end{equation*}
and 
\begin{equation*}
     \mathcal{K}_{a,b}=\left\{(u_n,A_n)_{n}\subset\mathcal{A}_{a,b} \left| \begin{array}{c} \rho_n\rightarrow \chi_{H^{-}_{a,b}} \text{ in } L^{p}(Q_{a,b}) \text{ for } p\geq 1 \\ B_n\rightarrow\frac{1}{\sqrt{2}}(1-\chi_{H^{-}_{a,b}}) \text{ in } L^{2}(Q_{a,b}) \end{array} \right. \right\}.
\end{equation*}
 For $0<\delta<a/2$, we also consider the sets
\begin{equation*}
\mathcal{A}^{\delta}_{a,b}=\left\{(u,A)\in\mathcal{A}_{a,b} \left| \begin{array}{c} (u,A)=(1,0) \text{ for } x_1<-\delta \\ (u,B)=(0,\frac{1}{\sqrt{2}}) \text{ for } x_1>\delta  \\ A_{1}=0 \end{array} \right.\right\}
\end{equation*}
and 
\begin{equation*}
    BC^{\delta}_{a,b}=\left\{(u_\n,A_\n)_\n\in\mathcal{K}_{a,b}\, |\, (u_n,A_n)\in \mathcal{A}^\delta_{a,b} \right\}.
\end{equation*}
We define the domain-dependent constants
 \begin{equation*}
        \sigma_{\kappa}(a,b)=\inf \left\{\liminf_{\n} \mathcal{E}_{\eps_n}(u_{\n},A_{\n},Q_{a,b}), \eps_n\rightarrow 0, (u_\n,A_\n)_\n\in\mathcal{K}_{a,b} \right\} 
    \end{equation*}
and
 \begin{equation*}
        \sigma_{\kappa,\delta}(a,b)=\inf \left\{\liminf_{n} \mathcal{E}_{\eps_n}(u_{\n},A_{\n},Q_{a,b}), \eps_n\rightarrow 0, (u_\n,A_\n)_\n\in BC_{a,b}^\delta \right\}.
    \end{equation*}
Notice that $\sigma_\kappa=\sigma_k
(1,1)$. We prove the following scaling property for $\sigma_\kappa$.  
\begin{proposition}

\label{prop:homogeneity}
        For any $a,b>0$, $N\in\N$ and $0<\delta<1/2$   
  \begin{equation}
  \label{e:homogeneity}
      \sigma_\kappa(a,b)=b\sigma_\kappa,
  \end{equation} 
\end{proposition}
and 
\begin{equation}\label{e:sigmabarest}
        \sigma_{\kappa,\delta}(1,N)\geq N\sigma_{\kappa,\delta}.
\end{equation}
       
\begin{proof}
   Let us start with \eqref{e:homogeneity}. Fix $L>0$ and consider a sequence $(u_\n,A_\n)_\n\in \mathcal{K}_{La,Lb} $ such that 
    \begin{equation*}
        \mathcal{E}_{\eps_n}(u_\n,A_\n,Q_{La,Lb})\rightarrow \sigma_{\kappa}(La,Lb).
    \end{equation*}
    Consider, for $S>0$, the change of variable $x=(L/S)y$, and define the sequence $\hat{\eps}_n=(S/L) \eps_{n}$, and the functions $\hat{u}_{n} (y)=u_\n(L/S y)$ and $\hat{A}_{\n} (y)=(S/L) A_\n(L/S y)$. We find that 
    \begin{equation*}
         \mathcal{E}_{\eps_n}(u_\n,A_\n,Q_{La,Lb})=\frac{L}{S}\mathcal{E}_{\hat{\eps}_n}(\hat{u}_\n,\hat{A}_\n,Q_{Sa,Sb}).
    \end{equation*} 
    This inequality implies, letting $\n\rightarrow \infty$,  \begin{equation*}
        \sigma_{\kappa}(La,Lb)\geq \frac{L}{S}\sigma_{\kappa}(Sa,Sb),
    \end{equation*}
    from which we deduce that, for any $L>0$,
    \begin{equation}
        \label{e:scaling}
        \sigma_{\kappa}(La,Lb)=L\sigma_{\kappa}(a,b).
    \end{equation} 
    The energy density is positive, thus, for  $(u_\n,A_\n)_\n\in \mathcal{K}_{a,b} $ and any $c,d>0$ with $Q_{c,d} \subset Q_{a,b}$, we have 
    \begin{equation*}
         E_{\eps_{\n}}(u_{\n},A_{\n},Q_{a,b})\geq E_{\eps_{\n}}(u_{\n},A_{\n},Q_{c,d}),
    \end{equation*}
    which implies 
    \begin{equation}
        \label{e:positiveenergy}
        \sigma_{\kappa}(a,b)\geq \sigma_{\kappa}(c,d).
    \end{equation}Moreover, for any $c,d>0$ and $(u_\n,A_\n)_\n\in \mathcal{K}_{a,c+d}$, by the translation invariance of the energy we have 
    \begin{equation*}
   E_{\eps_{\n}}(u_{\n},A_{\n},Q_{a,c+d})=E_{\eps_{\n}}(u_{\n},A_{\n},Q_{a,c})+E_{\eps_{\n}}(u_{\n},A_{\n},Q_{a,d}),
    \end{equation*}
    which implies 
    \begin{equation}
        \label{e:sumine}
        \sigma_\kappa(a,c+d)\geq \sigma_\kappa(a,c)+\sigma_\kappa(a,d).
    \end{equation}
    We now claim that 
    \begin{equation}
    \label{e:claimhom}
        \sigma_{\kappa}(1,b)=b \sigma_{\kappa}.
    \end{equation} Then, combining \eqref{e:claimhom} and \eqref{e:scaling}, we would get
    \begin{equation*}
        \sigma_{\kappa}(a,b)=a\sigma_{\kappa}\left(1,\frac{b}{a}\right)=b\sigma_{\kappa}
    \end{equation*}
     and the thesis would follow. Let us now prove \eqref{e:claimhom}.
    First, we show that \eqref{e:claimhom} holds for $b\in \N$. If $b\in\N$, from \eqref{e:sumine} we find immediately that $\sigma_\kappa(1,b)\geq b\sigma_\kappa$.
     For the opposite inequality, using \eqref{e:scaling} and \eqref{e:positiveenergy} we have 
    \begin{equation*}
        b\sigma_\kappa=\sigma_\kappa(b,b)\geq \sigma_\kappa(1,b).
    \end{equation*}   
    Next, we show that \eqref{e:claimhom} holds for $b=p/q$, with $p,q\in\N$.
    By \eqref{e:positiveenergy}, we have that $\sigma_{\kappa}(q,p)\geq \sigma_{\kappa}(1,p)$. Hence, using \eqref{e:scaling} we find 
    \begin{equation*}
        \sigma_{\kappa}(1,b)= \frac{1}{q}\sigma_{\kappa}(q,p)\geq \frac{1}{q}\sigma_{\kappa}(1,p)=b\sigma_{\kappa}.
    \end{equation*}
    Now, if $b>1$, we have using \eqref{e:scaling} and \eqref{e:positiveenergy}
    \begin{equation*}
        \sigma_{\kappa}(1,b)=b\sigma_{\kappa}\left(\frac{1}{b},1\right)\leq b \sigma_{\kappa}
    \end{equation*}
     and \eqref{e:claimhom} holds. 
     If $0< b <1$, by \eqref{e:sumine} we have
     \begin{equation}
     \label{e:finalcase}
         \sigma_{\kappa}(1,b+1)\geq \sigma_{\kappa}+\sigma_{\kappa}(1,b).
     \end{equation} Using \eqref{e:claimhom} for $b+1>1$, we have 
     \begin{equation}
         \label{e:finalcase2}
         \sigma_{\kappa}(1,b+1)= (b+1)\sigma_\kappa.
     \end{equation} Combining \eqref{e:finalcase} and \eqref{e:finalcase2} we can conclude that $b\sigma_{\kappa}\geq \sigma_{\kappa}(1,b)$, hence \eqref{e:claimhom} holds also for $b<1$. 
     
    Finally, the case $b\in\R^+$ follows by density and monotonicity of the function $b \mapsto \sigma_{\kappa}(1,b)$.
     
    The proof of \eqref{e:sigmabarest} follows from the same reasoning as \eqref{e:sumine}. For any $a,c,d>0$, by the translation invariance of the energy we have 
  \begin{equation*}
        \sigma_{\kappa,\delta}(a,c+d)\geq \sigma_{\kappa,\delta}(a,c)+\sigma_{\kappa,\delta}(a,d),
    \end{equation*}
    which implies \eqref{e:sigmabarest}.
\end{proof}
The goal now is to construct an appropriate minimizing sequence for $\sigma_\kappa$ such that we can compare it with a competitor for $\sigma_{\kappa,\delta}$. To do so, the following two lemmas are essential. Informally, we want to show that, given $\delta>0$ and $\eps_n\rightarrow 0$, we can ``cut" an appropriate competitor $(u_n,A_n)$ for $\sigma_\kappa$ in two points $t_1,t_2\in[-\delta,\delta]$ and ``glue" the prescribed boundary conditions of $\sigma_{\kappa,\delta}$. Such points need to be chosen such that $E_{\eps_{\n}}(u_n,A_n)$ is small on the segments $\{t_i\}\times (-1/2,1/2)$, $i=1,2$. This is similar to standard Modica-Mortola type arguments. However, in our setting extra care is necessary to handle the second term in \eqref{e:Aderformula}. In particular, we would like that on the segments $\{t_1\}\times (-1/2,1/2)$ the phase $\theta$ of $u$ is well defined. Since we are not able to ensure that this happens on a cube of height $1$, we perform a similar construction on a larger domain $Q_{1,N}$, with $N\in\N$ and $N>2$.
Using the scaling properties of $\sigma_\kappa$, this is enough to prove Proposition \ref{prop:sigmaexchange}.   

We start by constructing a minimizing sequence for which an appropriate Modica-Mortola type functional in $\rho$ vanishes away from the boundary and the interface. 
\begin{lemma}

\label{lem:striscia2}
    Let $0<\delta<1/2$, $\eps_n\rightarrow 0$ and $N\in \N$, $N>2$. Then, there exist a sequence of smooth functions $(u_{\n},A_{\n})_\n\in \mathcal{K}_{1,N} $ and values $ a_\n, b_\n, c_\n, d_\n \in \R$ such that  $ a_\n\in(-\delta,-3\delta/4)$, $ b_\n\in(-\delta/2,-\delta/4)$, $ c_\n\in(-N/2+\delta,-N/2+2\delta)$, $  d_\n\in (N/2-\delta,N/2-2\delta)$, $\mathcal{E}_{\eps_n}(u_{\n},A_{\n})\rightarrow \sigma_{\kappa}(1,N)$,
   \begin{equation}
   \label{e:striscia2}
        M_{\eps_{\n}}(\rho_\n,R_n)=\int_{R_n}{\eps_{\n}}\abs{\nabla \rho_{\n}}^{2}+\frac{1}{{\eps_{\n}}}W(\rho_{\n})\,dx\rightarrow 0,
    \end{equation}
    where $R_n=[ a_\n, b_\n]\times [ c_\n,  d_\n]$ and $W(\rho)=\frac{1}{2}\min(2\rho^{2},1)(1-\rho^{2})$, and
    \begin{equation}\label{zeroener}
     \mathcal{E}_{\eps_n}(u_{\n},A_{\n},R_n)\rightarrow 0
    \end{equation}

\end{lemma}
\begin{remark}
    The lemma follows from Lemma \ref{lem:mmestimate} and from the idea that no energy is present far from the interface. The main difficulty is the absence of periodic boundary conditions in this case, so we need to restrict the energy to a ``good'' rectangle $R$ to control the boundary term in \eqref{e:dwcontrol}.  The rectangle $R$ in principle depends on $\n$, but can be chosen such that it never degenerates for any $n\in \N$ and such that it never touches the boundary or the interface of $Q_{1,N}$, see Figure \ref{fig:rectangle}. This motivates the assumptions on the values $ a, b, c,  d$.
\end{remark}
\begin{proof}
   In this proof all implicit constants may depend on $N$ and $\kappa$. Let  $(u_{\n},A_{\n})_\n\in \mathcal{K}_{1,N} $ be such that $\mathcal{E}_{\eps_n}(u_{\n},A_{\n})\rightarrow \sigma_{\kappa}(1,N)$. By density of smooth functions in $\mathcal{A}(Q_{1,N})$, we can assume that $(u_{\n},A_{\n})_\n$ are smooth. Consider the set $(-\delta,-3\delta/4)\times (-N/2,N/2)$. Since $(u_{\n},A_{\n})$ is a minimizing sequence, clearly \begin{equation*}
       \mathcal{E}_{\eps_n}(u_{\n},A_{\n})\lesssim 1.
   \end{equation*} Hence, by the mean value theorem, it is possible to find $ a_\n \in (-\delta,-3\delta/4)$ such that
   \begin{equation*}
        \int_{-N/2}^{N/2}{\eps_{\n}}\abs{\nabla_{\alpha A_{\n}( a_\n,x_2)}u_{\n}( a_\n,x_2)}^2\,dx_2\lesssim \delta^{-1}.
    \end{equation*} 
    In a similar way it is possible to define $ b_\n, c_\n$ and $  d_\n$ satisfying the hypothesis above such that, letting 
    \begin{equation*}
        R_n=[ a_\n, b_\n]\times [ c_\n,  d_\n],
    \end{equation*} 
    we have
    \begin{equation}
        \label{e:boundaryconvergence}
        \int_{\partial R_n}{\eps_{\n}}\abs{\nabla_{\alpha A_{\n}}u_{\n}}^2\,d\mathcal{H}^1\lesssim \delta^{-1}.
    \end{equation}
        \begin{figure}[h]
    \centering
    \includegraphics[scale = 1.5]{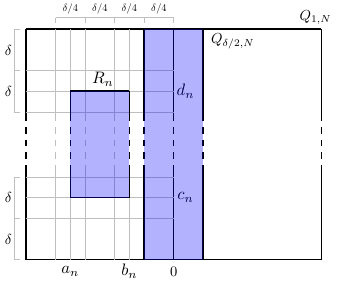} 
    \caption{Details of the sets $R_n$ and $Q_{\delta/2,N}$ inside $Q_{1,N}$.}
    \label{fig:rectangle}
\end{figure}
Let us prove \eqref{zeroener}. Since $R_n$ is bounded away from the interface uniformly in $\n$ (see Figure \ref{fig:rectangle}), by positivity of the energy
\begin{equation*}
         \mathcal{E}_{\eps_n}(u_{\n},A_{\n},Q_{1,N})\geq \mathcal{E}_{\eps_n}(u_{\n},A_{\n},R_n)+\mathcal{E}_{\eps_n}(u_{\n},A_{\n},Q_{\delta/2,N}).
     \end{equation*}
Passing to the limit we find 
     \begin{equation*}
         \sigma_{\kappa}(1,N)\geq \limsup_{\n} \mathcal{E}_{\eps_n}(u_{\n},A_{\n},R_n)+\sigma_{\kappa}(\delta/2,N).
     \end{equation*} Since by Proposition \ref{prop:homogeneity}, we know that $\sigma_{\kappa}(\delta/2,N)=\sigma_{\kappa}(1,N)$, this concludes the proof of \eqref{zeroener}. Now, letting $j_\n=j(u_\n,A_\n)$ (recall the definition in~\eqref{e:j}), by Lemma~\ref{lem:mmestimate} (and $\psi\lesssim 1$, recall \eqref{e:psi}) we have
\begin{equation*}
    M_{\eps_{\n}}(\rho_\n,R_n)\lesssim \mathcal{E}_{\eps_n}(u_\n,A_\n,R_n)+{\eps_{\n}}\int_{\partial R_n}\abs{j_\n} d\mathcal{H}^1 .
\end{equation*}
   Using that $\abs{j_\n}\leq \abs{\nabla_{\alpha A_\n} u_\n}$ since $\rho_\n\leq 1$, we find 
    \begin{equation*}
       {\eps_{\n}}\int_{\partial R_n}\abs{j_\n} d\mathcal{H}^1 \lesssim \eps_{\n}^\frac{1}{2}   \left(\int_{\partial R_n}{\eps_{\n}}\abs{\nabla_{\alpha A_{\n}}u_{\n}}^2\,d\mathcal{H}^1 \right)^{\frac{1}{2}},
    \end{equation*}
    from which the conclusion follows using \eqref{e:boundaryconvergence}.
    \end{proof}
Using Lemma \ref{lem:striscia2}, for any $\delta>0$, $\eps_n\rightarrow 0$ and $N\in \N$ with $N>2$ we now find an appropriate competitor $(u_n,A_n)$ for $\sigma_\kappa(1,N)$ and two values $t_1,t_2\in \R$ such that the energy is controlled  on the segments $\{t_i\}\times (-N/2+1,N/2-1)$, $i=1,2$ .
\begin{lemma}
\label{lem:t1t2}
Let $0<\delta<1/2$, ${\eps_{\n}}\rightarrow 0$ and $N\in \N$, $N>2$. Then, there exist a sequence of smooth functions $(u_{\n},A_{\n})_\n\in \mathcal{K}_{1,N} $ and values $t_1=t_1(\n)\in \R$, $t_2=t_2(\n)\in \R$ such that $t_{1}\in(-\delta,-\delta/4)$, $t_{2}\in(\delta/2,\delta)$, $\mathcal{E}_{\eps_n}(u_{\n},A_{\n})\rightarrow \sigma_{\kappa}(1,N)$ and 
\begin{enumerate}[label=(\roman*)]
    
    \item \label{it1} on $[-N/2+1,N/2-1]$, $\displaystyle \rho_{\n}(t_{1},\cdot)\rightarrow 1$ uniformly and  $\theta_n$ is well defined,
    \item \label{it3} $\displaystyle \int_{-N/2+1}^{N/2-1}{\eps_{\n}} \abs{\nabla \rho_{\n}(t_{1},x_{2})}^{2}dx_{2}\to 0 $,
    \item  \label{it5}$\displaystyle \int_{-N/2+1}^{N/2-1}{\eps_{\n}} \abs{\partial_2 \theta_{\n}(t_{1},x_{2})-\alpha A_{\n,2}(t_{1},x_{2})}^{2}dx_{2}\to 0$,
    \item \label{it8} $\displaystyle \int_{-N/2+1}^{N/2-1}{\eps_{\n}}\abs{\nabla_{\alpha A_{\n}}u_{\n}(t_{2},x_{2})}^{2} dx_{2}\to 0$,
    \item \label{it9} $\displaystyle \int_{-N/2+1}^{N/2-1}\abs{\rho_{\n}(t_{2},x_{2})}^p dx_{2} \rightarrow 0\qquad \forall p\geq 1$.
\end{enumerate}
\end{lemma}

\begin{proof}
    We start by taking the sequence $(u_{\n},A_{\n})_\n\in \mathcal{K}_{1,N} $ and the values $ a_\n, b_\n, c_\n,  d_\n \in \R$ constructed in Lemma \ref{lem:striscia2}. Arguing as for \eqref{zeroener}, we have
    \begin{equation*}
        \limsup_{n\to \infty} \mathcal{E}_{\eps_n}(u_{\n},A_{\n},(\delta/2,\delta)\times ( c_\n,  d_\n))=0.
    \end{equation*}
    Moreover, for all $p>1$
    \begin{equation*}
    \int_{Q_{1,N}\cap\{x_1\geq 0\}}\abs{\rho_n}^p\rightarrow 0,    
    \end{equation*}
    thus by the mean value theorem we can find $t_2\in (\delta/2,\delta)$ such that \textit{\ref{it8}} and \textit{\ref{it9}} hold. 
    Now, from Lemma \ref{lem:striscia2} we know that \eqref{e:striscia2} holds. Hence, again by the mean value theorem, we can find $t_1\in (a_\n, b_\n)\subset(-\delta,-\delta/4)$ such that \textit{\ref{it3}} holds, 
    \begin{align}
    \label{e:linetozero} &\int_{-N/2+1}^{N/2-1}{\eps_{\n}}\abs{\partial_2 \rho_{\n}(t_1,x_2)}^{2}+\frac{1}{\eps_\n} W(\rho_{\n}(t_1,x_2))\,dx_2\rightarrow 0, \\   
        \label{e:covariantt1} &\int_{-N/2+1}^{N/2-1}{\eps_{\n}}\abs{\nabla_{\alpha A_{\n}}u_{\n}(t_{1},x_{2})}^{2} dx_{2}\to 0
         \end{align}
     and 
     \begin{equation}
          \label{it6}  \int_{-N/2+1}^{N/2-1}\abs{\rho_{\n}(t_{1},x_{2})-1} dx_{2}\rightarrow 0.
     \end{equation}

    Concerning \textit{\ref{it1}}, take any $s,t\in(-N/2+1,N/2-1)$. We consider a real function $\eta$ such that $\eta'=2\sqrt{W}$. Using the usual Modica-Mortola trick and \eqref{e:linetozero} we find that 
    \begin{equation}
    \label{e:trick}
    \abs{\eta(\rho_{\n}(t_{1},t))-\eta(\rho_{\n}(t_{1},s))}\leq \int_{s}^{t}{\eps_{\n}}\abs{\partial_2 \rho_{\n}(t_1,x_2)}^{2}+\frac{1}{\eps_\n} W(\rho_{\n}(t_1,x_2))\,dx_2\rightarrow 0.    
    \end{equation} 
 Since $\eta$ is continuous and strictly increasing, \eqref{e:trick} and \eqref{it6} imply \textit{\ref{it1}}.
     
    Finally, by \textit{\ref{it1}}, the phase $\theta$ is well defined for $x_1=t_1$, so that by \eqref{e:covariantt1}
    \begin{equation}
         \label{it4} \int_{ -N/2+1}^{N/2-1}{\eps_{\n}} \rho_{\n}(t_{1},x_{2})^{2}\abs{\partial_2 \theta_{\n}(t_{1},x_{2})-\alpha A_{\n,2}(t_{1},x_{2})}^{2}dx_{2}\to 0.
    \end{equation} \textit{\ref{it1}} and \eqref{it4} together imply that \textit{\ref{it5}} holds.
\end{proof}
Finally, we can prove Proposition \ref{prop:sigmaexchange}.

\begin{proof}[Proof of Proposition \ref{prop:sigmaexchange}]
    Let $0<\delta<1/2$, $\eps_n\rightarrow 0$ and $N\in \N$ with $N>2$. Consider  $(u_{\n},A_{\n})_\n\in ~\mathcal{K}_{1,N} $, $t_1,t_2\in \R$ as given by Lemma \ref{lem:t1t2}.
    We claim that we can find a pair $(\hat{u}_\n,\hat{A}_\n)\in BC^{\delta}_{1,N-2}$ such that
    \begin{equation}
    \label{e:claim}
        \mathcal{E}_{\eps_n}(\hat{u}_\n,\hat{A}_\n,Q_{1,N-2})\leq \mathcal{E}_{\eps_n}(u_{\n},A_{\n},Q_{1,N})+o_\n(1).
    \end{equation}
    Then, by taking the liminf for $\n\rightarrow +\infty$ in \eqref{e:claim} we find $\sigma_{\kappa,\delta}(1,N-2)\leq \sigma_{\kappa}(1,N)$, so that by \eqref{e:homogeneity} and \eqref{e:sigmabarest} we have
    \begin{equation*}
        \sigma_{\kappa,\delta}\leq \frac{N}{N-2}\sigma_{\kappa}
    \end{equation*} and the result follows by taking the limit for $N\rightarrow +\infty $. Let us now construct a pair of functions that satisfy \eqref{e:claim}. To simplify notation, we drop the dependency on $\n$. By Lemma \ref{lem:t1t2} \textit{\ref{it1}}, for $n$ large enough on the segment $t_1\times({-N/2+1},N/2-1)$ we can write $u=\rho e^{i\theta}$ for a well-defined phase $\theta$. By fixing an appropriate gauge (recall \eqref{e:gauge}) we can assume that 
    $A=(0,A_2)$ and that on the segment
 $t_1\times({-N/2+1},N/2-1)$ we have $\theta = 0$ and thus $u=\rho$. Indeed, it suffices to take the change of gauge $\varphi(x_1,x_2)=\eta(x_1,x_2) -\theta(t_{1},x_{2})$ where $\eta$ is such that $\partial_1\eta=- \alpha A_{1}$ and $\eta(t_1,x_2)=0$ for $x_2\in (-N/2+1,N/2-1)$. We observe that the properties in Lemma \ref{lem:t1t2} are gauge invariant.
 
 We thus define for $x\in(-1/2,1/2)\times({-N/2+1},N/2-1)$
 \begin{equation*}
    \hat{u}(x_1, x_2)=
       \begin{cases}
       1& \text{ if } \ -\frac{1}{2}\leq x_{1}\leq t_{1}-\eps,  \\
       1-\frac{x_{1}-(t_{1}-\eps)}{\eps}(1 -\rho(t_{1},x_{2}))  & \text{ if } t_{1}-\eps \leq x_{1} \leq t_{1}, \\
       u(x_1, x_2)  &\text{ if } t_{1}\leq x_{1} \leq t_{2}, \\
       u(t_{2},x_{2})(1-\frac{x_{1}-t_{2}}{\eps})   &\text{ if } t_{2}\leq x_{1} \leq t_{2}+\eps, \\
       0  &\text{ if } t_{2}+\eps\leq x_{1} \leq \frac{1}{2} 
       \end{cases}
    \end{equation*}
 and $\hat{A}=(0,\hat{A}_2)$ where
\flushleft
 \begin{equation*}
    \hat{A}_2(x_1, x_2)=
       \begin{cases}
        0& \text{ if } \ -\frac{1}{2}\leq x_{1}\leq t_{1}-\eps,  \\
       A_{2}(t_{1},x_{2})\left(\frac{x_{1}-(t_{1}-\eps)}{\eps}\right)  & \text{ if } t_{1}-\eps \leq x_{1} \leq t_{1}, \\
       A_{2}(x_1, x_2)  &\text{ if } t_{1}\leq x_{1} \leq t_{2}, \\
        A_{2}(t_{2},x_{2})+\frac{1}{\sqrt{2}}(x_{1}-t_{2})  &\text{ if } x_{1} \geq t_{2}
       \end{cases}
    \end{equation*}

so that 
 \begin{equation*}
   \hat{B}(x_1, x_2)= \nabla \times \hat{A}(x_1, x_2)=
       \begin{cases}
       0 & \text{ if } \ -\frac{1}{2}\leq x_{1}\leq t_{1}-\eps,  \\
       \frac{1}{\eps}A_{2}(t_{1},x_{2}) & \text{ if } t_{1}-\eps \leq x_{1} \leq t_{1}, \\
       B(x_1, x_2)  &\text{ if } t_{1}\leq x_{1} \leq t_{2}, \\
       \frac{1}{\sqrt{2}}   &\text{ if } x_{1} \geq t_{2} .
       \end{cases}
    \end{equation*} 
By construction, we have $(\hat{u},\hat{A})\in BC^{\delta}_{1,N-2}$. Denoting \begin{equation*}
    f(x_{1})=\int_{{-N/2+1}}^{N/2-1}\left[ \eps\abs{\nabla_{\alpha \hat{A}}^\kappa \hat{u}}^{2}+\frac{1}{\eps}\left(\hat{B}-\frac{1}{\sqrt{2}}(1-\abs{\hat{u}}^{2})\right)^{2}\right](x_1,x_2)\,dx_{2},
\end{equation*} 
we have 
\begin{multline}
\label{e:splitting}
\mathcal{E}_\eps\left(\hat{u},\hat{A},\left(-\frac{1}{2},\frac{1}{2}\right)\times\left({-\frac{N}{2}+1},\frac{N}{2}-1\right)\right)=
\int_{-\frac{1}{2}}^{t_{1}-\eps}f(x_{1})dx_{1}+ \int_{t_{1}-\eps}^{t_{1}}f(x_{1})\,dx_{1}\\+\int_{t_{1}}^{t_{2}}f(x_{1})\,dx_{1}+ \int_{t_{2}}^{t_{2}+\eps}f(x_{1})\,dx_{1}+\int_{t_{2}+\eps}^{\frac{1}{2}} f(x_{1})\,dx_{1}.
\end{multline}
Let us examine individually each integral. We notice that the first and the last integral in \eqref{e:splitting} are equal to $0$. In the third integral $(\hat{u},\hat{A})=(u,A)$, hence we can estimate it by $\mathcal{E}_\eps(u,A,Q_{1,N})$. 
The rest of the proof is dedicated to estimating the second and the fourth integral in \eqref{e:splitting} using Lemma \ref{lem:t1t2}. Since $\abs{\nabla_{\alpha\hat{A}}^\kappa\hat{u}}^2\lesssim_\kappa\abs{\nabla_{\alpha\hat{A}}\hat{u}}^2$, we estimate the latter.


We analyze the second term in \eqref{e:splitting}. Recall that in $[t_1-\eps,t_1] \times (-1/2, 1/2)$ it holds $\hat{u}=\hat{\rho}$ so that  
\begin{equation*}
     \abs{\nabla_{\alpha\hat{A}}\hat{u}}^2=\abs{\nabla\hat{\rho}}^2+\alpha^2\hat{\rho}
     ^2\abs{\hat{A}_2}^2.
\end{equation*}
First, we have that for $x_1\in [t_1-\eps,t_1]$
\begin{equation*}
\hat{\rho}(x_1,x_2)= 1-\frac{x_{1}-(t_{1}-\eps)}{\eps}(1-\rho(t_{1},x_{2}))   
\end{equation*}
so that 
\begin{align}
\label{e:partial1}
    &\partial_1\hat{\rho}(x_1,x_2)=-\frac{1}{\eps}(1-\rho(t_{1},x_{2})),\\
 \label{e:partial2}
    &\partial_2\hat{\rho}(x_1,x_2)=\frac{x_{1}-(t_{1}-\eps)}{\eps}\partial_2\rho(t_{1},x_{2}).
\end{align}
Using \eqref{e:partial1}, \eqref{e:partial2} and Lemma \ref{lem:t1t2} \textit{\ref{it1}} and \textit{\ref{it3}} we find 
\begin{multline}
\label{e:t1grad}
    \int_{t_1-\eps}^{t_1}\int_{-N/2+1}^{N/2-1}\eps\abs{\nabla \hat{\rho}}^2\,dx\\\lesssim \int_{-N/2+1}^{N/2-1}\abs{1-\rho(t_1,x_2)}^2\,dx_2+ \eps^2\int_{-N/2+1}^{N/2-1}\abs{\partial_2\rho(t_1,x_2)}^2\,dx_2 = o_\eps(1).
\end{multline}
Furthermore, 
\begin{multline*}
   \eps \alpha^2 \int_{t_1-\eps}^{t_1}\int_{-N/2+1}^{N/2-1}\hat{\rho}
    ^2\abs{\hat{A}_2}^2\,dx\\ \lesssim \eps \alpha^2 \int_{t_1-\eps}^{t_1}\int_{-N/2+1}^{N/2-1}
    \abs{A_2(t_1,x_2)}^2\,dx =\eps^2 \alpha^2 \int_{-N/2+1}^{N/2-1}
    \abs{A_2(t_1,x_2)}^2\,dx_2.
\end{multline*}
Thus, by Lemma \ref{lem:t1t2} \textit{\ref{it5}} and our choice of gauge,
\begin{equation}
    \label{e:t1covar}
     \eps \alpha^2 \int_{t_1-\eps}^{t_1}\int_{-N/2+1}^{N/2-1}\hat{\rho}
    ^2\abs{\hat{A}_2}^2\,dx= o_\eps( \eps).
\end{equation}
Regarding the double-well part, using Young's inequality and $\hat{\rho}\leq 1$ we have 
\begin{multline}
\label{e:t1energ}
     \int_{t_1-\eps}^{t_1}\int_{-N/2+1}^{N/2-1}\frac{1}{\eps}\left(\hat{B}-\frac{1}{\sqrt{2}}(1-\hat{\rho}^{2})\right)^{2}\,dx \\\lesssim_\kappa   \int_{t_1-\eps}^{t_1}\int_{-N/2+1}^{N/2-1} \eps \alpha^2\abs{A_{2}(t_{1},x_{2})}^2+\frac{1}{    \eps}(1-\rho(t_1,x_2))^2\,dx \\ = \int_{-N/2+1}^{N/2-1} \eps^2 \alpha^2\abs{A_{2}(t_{1},x_{2})}^2+(1-\rho(t_1,x_2))^2\,dx_2 =o_\eps(1),
\end{multline}
where in the last step we used Lemma \ref{lem:t1t2} \textit{\ref{it1}} and \textit{\ref{it5}}. By \eqref{e:t1grad},\eqref{e:t1covar} and \eqref{e:t1energ} we can conclude that the second term in \eqref{e:splitting} is $o_\eps(1)$. We finally argue similarly for the fourth term in \eqref{e:splitting}.
In this case for  $x_1\in [t_2,t_2+\eps]$ and $x_2\in (-N/2+1,N/2-1)$, we have
\begin{align}
\label{e:partial3}
&\partial_1 \hat{u}(x_1,x_2)=-\frac{1}{\eps} u(t_{2},x_{2})\\
\label{e:partial4}
&\partial_2\hat{u}(x_1,x_2)=\partial_2 u(t_{2},x_{2})\left(1-\frac{x_{1}-t_{2}}{\eps}\right).
\end{align}
Lemma \ref{lem:t1t2} \textit{\ref{it9}} and \eqref{e:partial3} imply that 
\begin{equation}
\label{e:grad1estt2}
    \int_{t_2}^{t_2+\eps}\int_{-N/2+1}^{N/2-1} \eps \abs{\partial_1\hat{u}}^2\,dx =o_\eps(1).
\end{equation}
Now, \eqref{e:partial4} implies that for  $x_1\in [t_2,t_2+\eps]$ and $x_2\in (-N/2+1,N/2-1)$,
\begin{multline}
\label{e:covariantt2}
    \abs{\partial_2\hat{u}-i\alpha\hat{A}_2\hat{u}}^2(x_1,x_2)=\left(1-\frac{x_{1}-t_{2}}{\eps}\right)^2\abs{\partial_2u(t_2,x_2)-i\alpha\hat{A}_2(x_1,x_2)u(t_2,x_2)}^2\\\lesssim
    \abs{\partial_2u(t_2,x_2)-i\alpha A_{2}(t_2,x_2)u(t_2,x_2)}^2+ \alpha^2(x_1-t_2)^2\abs{u(t_2,x_2)}^2\\ \lesssim_{\kappa} \abs{\partial_2u(t_2,x_2)-i\alpha A_{ 2}(t_2,x_2)u(t_2,x_2)}^2+ \eps^{-2}\abs{u(t_2,x_2)}^2.
\end{multline}

Hence \eqref{e:covariantt2} and Lemma \ref{lem:t1t2} \textit{\ref{it8}} and \textit{\ref{it9}} imply that
\begin{equation}
\label{e:grad2estt2}
    \int_{t_2}^{t_2+\eps}\int_{-N/2+1}^{N/2-1} \eps \abs{\partial_2\hat{u}-i\alpha\hat{A}_2\hat{u}}^2\,dx =o_\eps(1).
\end{equation}
Moreover, notice that using Lemma \ref{lem:t1t2} \textit{\ref{it9}} we find
\begin{equation}
    \label{e:t2well}
     \int_{t_2}^{t_2+\eps}\int_{-N/2+1}^{N/2-1}\frac{1}{\eps}\left(\hat{B}-\frac{1}{\sqrt{2}}(1-\abs{\hat{u}}^{2})\right)^{2}\,dx=o_\eps(1).
\end{equation}
Combining \eqref{e:grad1estt2},\eqref{e:grad2estt2} and \eqref{e:t2well} we find that the fourth term in \eqref{e:splitting} is $o_\eps(1)$ as desired.

 \end{proof}

 \subsection{Characterization of $\sigma_\kappa$: horizontal periodic boundary conditions}
 \label{subsec:periodic}

  We now prove  that up to a small cost we can impose periodic boundary conditions for $u$ on the sides orthogonal to the transition line. In particular we want $u$ to coincide with a real valued one dimensional profile. To this extent, consider the profile (this choice is relatively arbitrary)
  \begin{equation*}
      U(t)=\begin{cases}
          1 &\text{ if }t\leq -1,\\
          -t & \text{ if } -1<t<0,\\
          0  & \text{ if }t\geq 0.
      \end{cases}
  \end{equation*}
  Now, given $\eps>0$, let us define $U_{\eps}:(-1/2,1/2)\rightarrow [0,1]$ as 
       \begin{equation}
       \label{e:u1d}
         U_{\eps}(t)=U\left(\frac{t}{\eps}\right).
     \end{equation}
     Notice that $U_{\eps}$ satisfies $U_{\eps}(t)=1$ for $t<-\eps$, $U_{\eps}(t)=0$ for $t>0$, 
  \begin{equation}
      \int_{-\frac{1}{2}}^\frac{1}{2}\eps \abs{ \dot{U}_{\eps}(t)}^2\,dt=1\label{e:1denergycontrol},
  \end{equation}
  and 
  \begin{equation}
  \label{e:1dconv}
      \int_{-\frac{1}{2}}^0\abs{U_{\eps}(t)-1}^2\,dt\leq \eps.
  \end{equation}

 \begin{proposition}
\label{prop:periodic}     
 
     Let $0<\delta<1/2$, then 
     \begin{equation}
     \label{e:sigmaineq}
         \sigma_{\kappa,\delta}^{\text{per}}\leq \sigma_{\kappa}+C_\kappa\delta,
     \end{equation}
     where 
     \begin{multline}
     \label{e:sigmaperdelta}
     \sigma_{\kappa,\delta}^{\text{per}}=\inf\left\{\liminf_\n \mathcal{E}_{\eps_n}(u_\n,A_\n)|\eps_n\rightarrow 0, (u_\n,A_\n)_\n\in BC^{\delta}, \right. 
     \\ \left. u_\n(x_1,x_2)=U_{{\eps_{\n}}}(x_1) \text{ for }x_2=\pm 1/2 \right\}.
     \end{multline}
\end{proposition}
 \begin{proof}
    In this proof all implicit constants depend on $\kappa$. Let $\delta\in (0,1/4)$ and $\eps_n\rightarrow 0$. By \eqref{e:sigmabar}, we can find a sequence $(u_\n,A_\n)_\n\in BC^\delta$ such that $E_{\eps_{\n}}(u_\n,A_\n)\rightarrow\sigma_{\kappa,\delta}$. We prove that there exists a sequence $(\hat{u}_{\n},\hat{A}_\n)_\n\in BC^{2\delta}$ with $ \hat{u}_\n=U_{{\eps_{\n}}} \text{ for }x_2=\pm 1/2$ and such that
     \begin{equation}
     \label{e:goal}
         \mathcal{E}_{\eps_n}(\hat{u}_{\n},\hat{A}_\n)\leq \mathcal{E}_{\eps_n}(u_\n,A_{\n})+C\left(\delta +\frac{\eps_n}{\delta}\right),
     \end{equation}
        which, combined with Proposition \ref{prop:sigmaexchange}, implies \eqref{e:sigmaineq}. We focus on the top edge $x_2=1/2$, the same reasoning applies for the bottom edge. Since $\mathcal{E}_{\eps_n}(u_{\n},A_{\n})\lesssim 1$, by the mean value theorem we can find $1/2-\delta<\overline{t}<1/2$ such that \begin{equation}
      \label{e:deltaenergyest}
          \int_{-\frac{1}{2}}^\frac{1}{2}{\eps_{\n}}  \abs{\partial_1 u_\n(x_1,\overline{t})}^2\,dx_1\lesssim\frac{1}{\delta}.
      \end{equation}  We can  
     further assume that ${\eps_{\n}}<\delta$ and that $\bar{t}+{\eps_{\n}}<1/2$.
     We define the function $\hat{u}_\n$ as  
      \begin{equation*}
     \hat{u}_\n(x_1,x_2)=
         \begin{cases}
             u_\n(x_1,x_2) & \text{ if }  x_2\leq\overline{t},\\
              u_\n(x_1,\overline{t})\frac{\overline{t}+{\eps_{\n}}-x_2}{\eps_\n} +U_{\eps_{\n}}(x_1+\delta)\frac{x_2-\overline{t}}{{\eps_{\n}}}&\text{ if } \overline{t}\leq x_2\leq\overline{t}+ {\eps_{\n}},\\
              U_{\eps_{\n}}(x_1+\delta)\frac{\overline{t}+2{\eps_{\n}}-x_2}{{\eps_{\n}}} +U_{\eps_{\n}}(x_1)\frac{x_2-\overline{t}-{\eps_{\n}}}{{\eps_{\n}}}&\text{ if } \overline{t}+{\eps_{\n}}\leq x_2\leq\overline{t}+ 2{\eps_{\n}},\\
            U_{{\eps_{\n}}}(x_1) & \text{ if }  x_2\geq \bar{t}+2{\eps_{\n}}.\\

         \end{cases}
     \end{equation*}
     Similarly, we define the vector field $\hat{A}_\n=(0,\hat{A}_{\n,2})$ where
     \begin{equation*}
         \hat{A}_{\n,2}(x_1,x_2)=
         \begin{cases}
             A_{\n,2}(x_1,x_2) &\text{ if } x_2\leq\overline{t},\\
             \max \{\frac{1}{\sqrt{2}}(x_1-\delta),0\} &\text{ if } \overline{t}\leq x_2\leq \overline{t}+2{\eps_{\n}},\\
             \max \{\frac{1}{\sqrt{2}}x_1,0\} &\text{ if } x_2\geq \overline{t}+2{\eps_{\n}}.
         \end{cases}
     \end{equation*}
     Notice that $\hat{A}_{\n,2}$ has a discontinuity, but still $\nabla \times \hat{A}_{\n}=\partial_1\hat{A}_{\n,2}\in L^2(Q_1)$. As a result, it holds  $(\hat{u}_\n,\hat{A}_\n)\in BC^{2\delta}$. We now prove estimate \eqref{e:goal}. Letting $F_1=\{x_2\leq \overline{t}\}$, $F_2=\{\overline{t}\leq x_2\leq \overline{t}+{\eps_{\n}}\}$, $F_3=\{\overline{t}+{\eps_{\n}}\leq x_2\leq \overline{t}+2{\eps_{\n}}\}$ and $F_4=\{x_2\geq \overline{t}+2{\eps_{\n}}\}$, we have 
     \begin{equation}\label{e:split}
         \mathcal{E}_{\eps_n}(\hat{u}_{\n},\hat{A}_\n,Q_1)=\mathcal{E}_{\eps_n}(\hat{u}_{\n},\hat{A}_\n,F_1)+\mathcal{E}_{\eps_n}(\hat{u}_{\n},\hat{A}_\n,F_2)+\mathcal{E}_{\eps_n}(\hat{u}_{\n},\hat{A}_\n,F_3)+\mathcal{E}_{\eps_n}(\hat{u}_{\n},\hat{A}_\n,F_4). 
         \end{equation}
     We estimate each term in \eqref{e:split} separately. On $F_1$ it holds $(\hat{u}_{\n},\hat{A}_\n)=(u_{\n},A_\n)$, hence 
     \begin{equation}
     \label{e:bulkcontrol}
         \mathcal{E}_{\eps_n}(\hat{u}_{\n},\hat{A}_\n,F_1)\leq \mathcal{E}_{\eps_n}(u_{\n},A_\n,Q_1). 
     \end{equation}
    On $F_2$, using the fact that $\abs{\nabla_{\hat{A}}^\kappa \hat{u}}^2 \lesssim \abs{\nabla_{\hat{A} }\hat{u}}^2$, we have
     \begin{equation}
     \label{e:transition}
          \mathcal{E}_{\eps_n}(\hat{u}_{\n},\hat{A}_\n,F_2)\lesssim \int_{\overline{t}}^{\overline{t}+{\eps_{\n}}}\int_{-\delta-{\eps_{\n}}}^{\delta}{\eps_{\n}}\left(\abs{\partial_1 \hat u_{\n}}^{2}+\abs{\partial_2 \hat{u}_{\n}}^{2}\right)+\frac{1}{2{\eps_{\n}}}(1-\abs{\hat{u}_\n}^{2})^2\,dx_1dx_2.
     \end{equation}
    By definition of $\hat{u}_\n$, on $F_2$ we find \begin{align}
     &\partial_1\hat{u}_\n(x_1,x_2)=\partial_1u_{n}(x_1,\overline{t})\frac{\overline{t}+{\eps_{\n}}-x_2}{{\eps_{\n}}} +\dot{U}_{\eps_{\n}}(x_1+\delta)\frac{x_2-\overline{t}}{{\eps_{\n}}}\label{e:partial1u},\\
     &\partial_2\hat{u}_\n(x_1,x_2)=\frac{1}{{\eps_{\n}}}(U_{\eps_n}(x_1+\delta)- u_\n(x_1,\overline{t}))\label{e:partial2u}.
     \end{align}
     Equations \eqref{e:partial1u}, \eqref{e:1denergycontrol} and \eqref{e:deltaenergyest} imply that 
     \begin{multline}
     \label{e:partial1control}
          \int_{\overline{t}}^{\overline{t}+{\eps_{\n}}}\int_{-\delta-{\eps_{\n}}}^{\delta}{\eps_{\n}}\abs{\partial_1 \hat u_{\n}}^{2}dx_1dx_2\\ \lesssim   \int_{\overline{t}}^{\overline{t}+{\eps_{\n}}}\int_{-\delta-{\eps_{\n}}}^{\delta}{\eps_{\n}}\left(\abs{\partial_1 u_{\n}(x_1,\overline{t})}^{2}+\abs{\dot{U}_{\eps_{\n}}(x_1)}^2\right)dx_1dx_2\lesssim\frac{{\eps_{\n}}}{\delta}+\eps_n\lesssim \frac{\eps_n}{\delta}.
     \end{multline}
     Furthermore, \eqref{e:partial2u} implies that 
     \begin{multline}
     \label{e:partial2control}
          \int_{\overline{t}}^{\overline{t}+{\eps_{\n}}}\int_{-\delta-{\eps_{\n}}}^{\delta}{\eps_{\n}}\abs{\partial_2 \hat u_{\n}}^{2}dx_1dx_2 \\ \lesssim   \int_{\overline{t}}^{\overline{t}+{\eps_{\n}}}\int_{-\delta-{\eps_{\n}}}^{\delta}\frac{1}{{\eps_{\n}}}(\abs{u_{\n}(x_1,\overline{t})}^{2}+\abs{U_{\eps_{\n}}(x_1+\delta)}^2)dx_1dx_2\lesssim \delta.
     \end{multline}
     Moreover, since $\abs{\hat{u}_n}\leq 1$, 
     \begin{equation}
         \label{e:wellcontrol}
         \int_{\overline{t}}^{\overline{t}+{\eps_{\n}}}\int_{-\delta-{\eps_{\n}}}^{\delta}\frac{1}{2{\eps_{\n}}}(1-\abs{\hat{u}_\n}^{2})^2\,dx_1dx_2\lesssim \delta.
     \end{equation}
     Plugging \eqref{e:partial1control},\eqref{e:partial2control} and \eqref{e:wellcontrol} in \eqref{e:transition} we get 
     \begin{equation}
         \label{e:transition control}
         \mathcal{E}_{\eps_n}(\hat{u}_{\n},\hat{A}_\n,F_2)\lesssim \delta +\frac{\eps_n}{\delta}.
     \end{equation}
     Let us now consider $F_3$. With a similar computation as before 
          \begin{multline*}
          \int_{\overline{t}+{\eps_{\n}}}^{\overline{t}+2{\eps_{\n}}}\int_{-\delta-{\eps_{\n}}}^{\delta}{\eps_{\n}}\abs{\partial_1 \hat u_{\n}}^{2}dx_1dx_2 \\ \lesssim   \int_{\overline{t}+{\eps_{\n}}}^{\overline{t}+2{\eps_{\n}}}\int_{-\delta-{\eps_{\n}}}^{\delta}{\eps_{\n}}\left(\abs{\dot{U}_{\eps_{\n}}(x_1+\delta)}^2+\abs{\dot{U}_{\eps_{\n}}(x_1)}^2\right)dx_1dx_2\lesssim {\eps_{\n}}
     \end{multline*}
    and
    \begin{equation*}
    \int_{\overline{t}+{\eps_{\n}}}^{\overline{t}+2{\eps_{\n}}}\int_{-\delta-{\eps_{\n}}}^{\delta}{\eps_{\n}}\abs{\partial_2 \hat u_{\n}}^{2}+\frac{1}{2{\eps_{\n}}}(1-\abs{\hat{u}_\n}^{2})^2\,dx_1dx_2\lesssim \delta
    \end{equation*}
     so that 
     \begin{equation}
         \label{e:shiftcontrol}
         \mathcal{E}_{\eps_n}(\hat{u}_{\n},\hat{A}_\n,F_3)\lesssim  \delta+\eps_n.
     \end{equation}
      Finally we address $F_4$. By \eqref{e:u1d}, \eqref{e:1denergycontrol}, \eqref{e:1dconv} and the definitions of $\hat{u}_\n$ and $\hat{A}_\n$ we have
     \begin{equation}
     \label{e:boundarycontrol}
         \mathcal{E}_{\eps_n}(\hat{u}_{\n},\hat{A}_\n,F_4)=\int_{\overline{t}+2{\eps_{\n}}}^\frac{1}{2}\int_{-{\eps_{\n}}}^0{\eps_{\n}}\abs{\dot{U}_{\eps_{\n}}}^{2}+\frac{1}{2{\eps_{\n}}}\left(1-\abs{U_{\eps_{\n}}}^{2}\right)^{2}\,dx_1dx_2\lesssim \delta.
     \end{equation}
     Combining \eqref{e:bulkcontrol}, \eqref{e:boundarycontrol}, \eqref{e:transition control} and \eqref{e:shiftcontrol} gives \eqref{e:goal}.
    \end{proof}

 \begin{remark}
 \label{rem:1D}
 Let $\rho,A:\R\rightarrow \R$ and consider for $a<b\in \R\cup\{\pm \infty\}$
 \begin{equation*}
     \mathcal{E}^{1d}(\rho,A,(a,b))=\int_a^b\abs{\dot{\rho}}^2+\kappa^{-2}\rho^2A^2+\left(\dot{A}-\frac{1}{\sqrt{2}}(1-\rho^2)\right)^2\,dt.
 \end{equation*}
     In \cite{chapman2000,chapman_91}, the following one dimensional phase transition problem for type-I superconductors is studied:
     \begin{equation*}
         \sigma^{1d}_\kappa=\inf\left\{E^{1d}(\rho,A,(-\infty,+\infty)),(\rho,A)(-\infty)=(1,0),\, (\rho,\dot{A})(+\infty)=\left(0,\frac{1}{\sqrt{2}}\right) \right\}.
     \end{equation*}
     Notice that if we consider $\Tilde{u}(x_1,x_2)=\rho(x_1)$ and $\Tilde{A}(x_1,x_2)=[0,A(x_1)]$, then 
     \begin{equation*}
\mathcal{E}_1(\Tilde{u},\Tilde{A},Q_1)=\mathcal{E}^{1d}\left(\rho,A,\left(-\frac{1}{2},\frac{1}{2}\right)\right),
     \end{equation*}
     which implies $\sigma_\kappa\leq \sigma_\kappa^{1d}$. We conjecture that $\sigma_\kappa= \sigma^{1d}_\kappa$, which would thus imply that the transition in our problem is one dimensional. Howevere due to the vectorial nature of our problem, this seems to be a difficult question.

 \end{remark}
 \section{$\Gamma$-limsup inequality}
 \label{sec:limsup}
 
 The goal of this section is to prove the following theorem (Theorem \ref{th:main}-\ref{i:limsup}).
 \begin{theorem}
 \label{th:limsup}
 
     Let $\eps_n\rightarrow 0$ such that $\eps_n^{-2}\kappa^{-2}b_{\text{ext}}\in 2\pi\Z$. Let $\rho \in BV(\tor,\{0,1\})$ be such that $\int_{\tor} (1-\rho^2)/\sqrt{2}\,dx=b_{\text{ext}}/\kappa\in(0,1)$. Then, there exists a sequence $(u_{n},A_{n})\in \mathcal{A}_{\eps_n}$  such that
    $\abs{u_{n}}=\rho_{n}\rightarrow \rho$ in $L^{p}(\tor)$ for all $p\geq 1$, $B_n=\nabla\times A_{n}\rightarrow (1-\rho^2)/\sqrt{2}$ in $L^{2}(\tor)$, $\int_{\tor} B_n\,dx=b_{\text{ext}}/\kappa$  and 
    \begin{equation}
    \label{e:gammalimsup}
        \limsup_n \mathcal{E}_{\eps_n}(u_n,A_n)\leq \sigma_\kappa P(\{\rho=0\},\tor). 
    \end{equation}
 \end{theorem}
A function $\rho\in BV(\tor,\{0,1\})$ is characterized by the set of finite perimeter $E=\{\rho=0\}$. To understand the main ideas behind the proof, we first start by the easier case where $E$ is a simply connected closed polyhedral set. We then discuss the generalization to the case where $E$ is a connected closed polyhedral set, then finally we examine the construction on a generic closed polyhedral set. The result for a set of finite perimeter follows by a standard density argument. In all these arguments the fundamental building block for the construction is given by the following proposition, which follows immediately from Proposition \ref{prop:sigmaexchange}, \eqref{e:sigmaineq} and a smoothing argument.
\begin{proposition}
\label{prop:buildingblock}
  Let $\xi>0$. Then, there exist $0<\eps_0<\delta_0<1/2$ and  $(u_0,A_0)\in \mathcal{A}^{\delta_0}$ such that $u_0=U_{\eps_0}$ for $x_{2}=\pm 1/2$ and 
    \begin{equation*}
        \mathcal{E}_{\eps_0}(u_0,A_0,Q_1)\leq \sigma_\kappa+\xi.
    \end{equation*}
We write $B_0=\nabla\times A_0$ and $\rho_0=\abs{u_0}$. Moreover, $u_0$ is continuous and there exists a function $\theta_0\in H^1_{\text{loc}}(\{\rho_0>0\})$ such that in the set $\{\rho_0>0\}$ we have $u_0=\rho_0e^{i\theta_0}$. In particular, on $\{x_1=-1/2\}$ and on $\{x_1\leq 0\}\cap \{x_2=\pm 1/2\}$ it holds $\rho_0>0$ and $\theta_0=0$. Furthermore,
    \begin{equation}
    \label{e:b0int}
        \frac{1}{4\sqrt{2}}<\int_{Q_1}B_0\,dx\leq \frac{3}{4\sqrt{2}}.
    \end{equation}
\end{proposition}
\begin{remark}
    Condition \eqref{e:b0int} comes from the convergence constraint on the magnetic field in \eqref{e:K}.
\end{remark}
\begin{proof}[Proof of Theorem \ref{th:limsup}]Let $E=\{\rho=0\}$. In this proof all implicit constants may  depend on $E$. By \cite[Remark 13.13]{maggi}, closed polyhedral sets are dense in finite perimeter sets with respect to the $L^1$ topology. Hence, as usual in $\Gamma$-convergence results (see \cite{braides,dalmaso,dalmasoleoni}), by a diagonal argument we may restrict the construction of the recovery sequence to closed polyhedral sets.

\quad \emph{\textbf{Case 1:  $E$ simply connected closed polyhedral set.} }
Let $E\subset \tor$ be a closed simply connected polyhedral set with $N$ edges. That is $\partial E=\cup_{i=1}^Ne_i$ where the sets $e_i$ are segments intersecting at their endpoints.  We denote, for $i=1\dots N$, by $L_i$ the length of $e_i$ and by $\nu_i$ the inward normal vector of $E$ on $e_i$. Let us also call $c_i^\pm$ the two endpoints of $e_i$.  \\
We proceed as follows: first, we define $\rho_n$ and $B_n$ that satisfy the required convergence assumptions and the mass constraint and then, we construct an appropriate couple $(u_n,A_n)\in \mathcal{A}_{\eps_n}$ that satisfies the energy bound. \\
Let $\xi>0$. Consider $(u_0,A_0)$, $\eps_0$, and $\delta_0$ as given by Proposition \ref{prop:buildingblock}. The idea is to partially cover each edge of $E$ with  mutually disjoint rectangular sets $\Omega_i$ of width $\eps_n/\eps_0$ and  length almost equal to $ L_i$. Moreover, if $e_i\subset \partial E_k$, we shift the sets $\Omega_i$ in the direction given by $\nu_i$ by a factor $\zeta=\zeta(\eps_n)\rightarrow 0$ to be chosen later. To do so, for any set $F\subset \tor$ consider the signed distance function
\begin{equation*}
    sd_F(x)=\begin{cases}
        d(x,\partial F) &\text{ if } x\in F,\\
        -d(x,\partial F) &\text{ if } x\in F^c.
    \end{cases}
\end{equation*} 
Let $\zeta \in (-\eps_n/(2\eps_0),\eps_n/(2\eps_0))$ to be chosen later and define for any $F\subset \tor$
\begin{equation}
\label{e:epsneigh}
     T(\eps_n,\zeta,F)=\left\{x\in Q_1,\, -\frac{\eps_n}{2\eps_0}+\zeta <  sd_F(x 
     )< \frac{\eps_n}{2\eps_0}+\zeta\right\}.
\end{equation} 
Notice that $T(\eps_n,0, E)$ is just an $ (\eps_n/\eps_0)$-neighborhood of $E$.
For $i=1\dots N$, we set
\begin{equation*}
    \Omega_i(\eps_n,\zeta)=\left\{x\in T(\eps_n,\zeta,E),\, d(x,e_i+\zeta \nu_i)< \frac{\eps_n}{2\eps_0},\,\abs{P_{e_i}x-c_i^+}> l_i^+,\abs{P_{e_i}x-c_i^-}> l_i^- \right\},
\end{equation*}
where $P_{e_i}$ is the projection on the line containing $e_i$ and $l_{i}^\pm=l_{i}^\pm(\eps_n)\rightarrow 0$ are such that $\Omega_i(2\eps_n,0)\cap \Omega_j(2\eps_n,0)=\emptyset$. This last condition implies that for any $\zeta\in(-\eps_n/(2\eps_0),\eps_n/(2\eps_0))$ it holds $\Omega_i(\eps_n,\zeta)\cap \Omega_j(\eps_n,\zeta)=\emptyset$. Moreover, it is possible to assume that $L_i-l_{i}^--l_{i}^+=N_i\eps_n/\eps_0$, with $N_i=N_i(\eps_n)\in\N$, so that each $\Omega_i(\eps_n,\zeta)$ consists of $N_i$ squares of side $\eps_n/\eps_0$, see Figure \ref{fig:angle}. For $j=1\dots N_i(\eps_n)$, let us call $x_{i,j}$ the center of the square $j$ constituting $\Omega_i$. For simplicity, we call $T=T(\eps_n,\zeta, E)$ and   $\Omega_i=\Omega_i(\eps_n,\zeta)$. Notice that, by construction, 
\begin{equation}
\label{e:smallangles}
    \abs{ T\backslash\bigcup_i \Omega_i}\lesssim \eps_n^2.
\end{equation}
We can now define $\rho_n$ and $B_n$ by ``gluing" an appropriate rescaling of the building blocks $(u_0,A_0)$ on the sets $\Omega_i$. 
Recalling \eqref{e:u1d}, let 
\begin{equation}
\label{e:rho}\rho_n(x)=
    \begin{cases}
        0 &\text{in}  \ E\backslash  T,\\
        \rho_0(\frac{\eps_0}{\eps_n}R_{\nu_i}(x-x_{i,j})) &\text{in} \ Q_{\frac{\eps_n}{\eps_0}}(x_{i,j}),\\
        U_{\frac{\eps_n}{\eps_0}}(sd_E(x)-\zeta) &\text{in} \  T\backslash \bigcup_i \Omega_i,\\
        1 &\text{elsewhere}
    \end{cases}  
\end{equation}
and 
\begin{equation}
\label{e:b}B_n(x)=
    \begin{cases}
        \frac{1}{\sqrt{2}} &\text{in}  \ E\backslash  T,\\
        B_0(\frac{\eps_0}{\eps_n}R_{\nu_i}(x-x_{i,j})) &\text{in} \ Q_{\frac{\eps_n}{\eps_0}}(x_{i,j}),\\
        0 &\text{elsewhere}.
    \end{cases}  
\end{equation}
The definition of $\rho_n$ on $ T\backslash \bigcup_i \Omega_i$ is just a rescaled Lipschitz extension of $U_{\eps_0}$, which is the horizontal boundary condition of $u_0$ (compare with \eqref{e:u1d}). Hence, $\rho_n$ is continuous. 
    \begin{figure}[h]
    \centering
    \includegraphics[scale = 0.9]{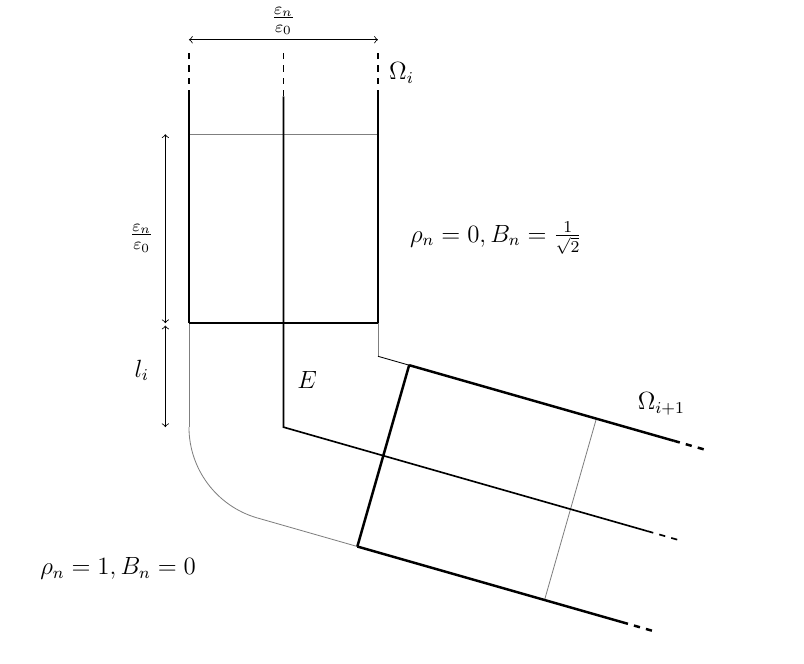} 
    \caption{Sketch of the construction near a corner of $\partial E$. For simplicity in this picture $\zeta=0$.}
    \label{fig:angle}
\end{figure}
A computation shows that $\rho_n$ and $B_n$ satisfy 
\begin{equation*}
    \int_{\tor}\abs{\rho_n-\rho}+\abs{B_n-\frac{1}{\sqrt{2}}(1-\rho)}^2\,dx\lesssim \abs{T}\lesssim \eps_n, 
\end{equation*}
which implies $\rho_n\rightarrow \rho\in L^1(\tor)$ (hence in $L^{p}(\tor)$ for all $1\leq p<\infty$ since, by construction, $\rho_n\leq 1$) and $B_n\rightarrow(1-\rho^2)/\sqrt{2}\in L^2(\tor)$, as required. \\
Regarding the flux constraint of $B_n$, we have
\begin{equation*}
    \int_{\tor}B_n\,dx=\int_{E\cup T}B_n\,dx=\frac{1}{\sqrt{2}}\abs{E\backslash T}+\frac{\eps_n^2}{\eps_0^2}\sum_i{N_i}\int_{Q_1}B_0\,dx=h(\zeta).
\end{equation*}
By construction we have that \begin{equation*}
    \abs{E\backslash T}= \abs{E}-P(E)\left(\zeta+\frac{\eps_n}{2\eps_0}\right)+o(\eps_n)
\end{equation*}
and $\sum_i{N_i}\eps_n/\eps_0\rightarrow P(E)$. Thus, 
\begin{equation}
\label{e:hest}
    h(\zeta)=\frac{1}{\sqrt{2}}\abs{E}-P(E)\left[\frac{1}{\sqrt{2}}\left(\zeta+\frac{\eps_n}{2\eps_0}\right)-\frac{\eps_n}{\eps_0}\int_{Q_1}B_0\,dx\right]+o(\eps_n).
\end{equation}
By \eqref{e:b0int}, for $n$ large enough  $ h(-\eps_n/(2\eps_0))>\abs{E}/\sqrt{2}$ and $ h(\eps_n/(2\eps_0))<\abs{E}/\sqrt{2}$. Since $h\in C^0$,  we can choose $\zeta \in (-\eps_n/(2\eps_0),\eps_n/(2\eps_0))$ such that  
\begin{equation}
\label{e:massconstr}
    \int_{\tor}B_n\,dx=\frac{1}{\sqrt{2}}\abs{E}=\frac{b_{ext}}{\kappa}.
\end{equation}  \\ 
Let us now construct  $(u_n,A_n)\in \mathcal{A}_{\eps_n}$ that satisfies the energy bound.\\
Since we have defined $B_n$ everywhere in $\tor$, we can find a global vector potential $A_n$ such that $\nabla \times A_n=B_n$. For instance we can take $A_n=(0,A_{n,2})$ with 
\begin{equation}
\label{e:aglobal}
    A_{n,2}(x_1,x_2)=\int_{-\frac{1}{2}}^{x_1}B_n(t,x_2)dt.
\end{equation}
We now construct an appropriate phase $\theta_n$ such that the non local term $ \nabla \theta_n-\alpha A_n$ is controlled. To this extent, we introduce an auxiliary vector field by setting for $x\in E^c\cup T$
\begin{equation}
\label{e:a_0}
A^0_n(x)=
    \begin{cases}
        \frac{\eps_n}{\eps_0}R_{\nu_i}^{-1}A_0(\frac{\eps_0}{\eps_n}R_{\nu_i}(x-x_{i,j})) &\text{in} \ Q_{\frac{\eps_n}{\eps_0}}(x_{i,j}),\\
        0 &\text{elsewhere}.
    \end{cases}  
\end{equation}
Notice that $A^0_n$ has discontinuities but still $\nabla\times A^0_n\in L^2(E^c\cup T)$, since the segments of discontinuity are perpendicular to the direction of $A^0_n$ (recall the definition of $\nabla\times$). We also define an auxiliary phase function.  For $x\in  \{\rho_n>0\}$, let 
\begin{equation}
\label{e:thetazero}
\theta^0_n(x)=
    \begin{cases}
        \theta_0(\frac{\eps_0}{\eps_n}R_{\nu_i}(x-x_{i,j})) &\text{in} \ Q_{\frac{\eps_n}{\eps_0}}(x_{i,j}),\\
        0 &\text{elsewhere}.
    \end{cases} 
\end{equation}
The function  $\theta^0_n$ is continuous in $\{\rho_n>0\}$ since $u_0=U_{\eps_0}$ on the horizontal boundaries of $Q_1$ and $u=1$ on the left boundary of $Q_1$, hence on $\partial Q_1\cap \{\rho_0>0\}$ $u_0$ is real and $\theta_0=0$. We remark that up to this point we have not used the hypothesis that $E$ is simply connected. \\
Since $E$ is closed and simply connected, the set $E^c\cup T$ is open and connected. Let $x_0\in E^c\cup T$ and for $x\in  \{\rho_n>0\}$ let $\Gamma_{x,x_0}\subset E^c\cup T$ be any Lipschitz curve connecting $x$ to $x_0$. We define for $x\in  \{\rho_n>0\}$
\begin{equation}
\label{e:theta}
    \theta_n(x)= \theta^0_n(x)+
    \alpha\int_{\Gamma_{x,x_{0}}}(A_n-A^0_n)\cdot dl.
\end{equation}
Let us prove that $\theta_n$ is well defined up to $2\pi \Z$. Let $\Gamma$ be a simple closed Lipschitz curve in $E^c\cup T$ enclosing $E\backslash  T$. Since in $E^c\cup T$ the term $A_n-A_n^0$ is curl free, by Stokes' theorem its path integral is invariant for homotopies of $\Gamma$ and in particular
\begin{equation*}
   \int_\Gamma (A_n-A^0_n)\cdot\,dl = \int_{\partial(E\cup  T)} (A_n-A^0_n)\cdot\,dl.
\end{equation*} 
Moreover, since  $A_n^0=0$ on $\partial(E\cup T)$, we have that 
\begin{equation}
\label{e:quantization}
    \alpha\int_{\partial(E\cup  T)} (A_n-A^0_n)\cdot\,dl =\alpha\int_{\partial(E\cup  T)} A_n\cdot\,dl =\alpha\int_{E\cup  T}B_n\,dx=\alpha \frac{b_{\text{ext}}}{\kappa}=\eps_n^{-2}\kappa^{-2}b_{\text{ext}}.
\end{equation}
Since by assumption $\eps_n^{-2}\kappa^{-2}b_{\text{ext}}\in 2\pi\Z$, \eqref{e:quantization} implies that $\theta_n(x)$ is defined up to a multiple of $2\pi$, and thus we can uniquely define  $u_n=\rho_ne^{i\theta_n}$.
Moreover, we remark that on $\{\rho_n>0\}$ it holds  
\begin{equation}
\label{e:thetacondition}
    \nabla \theta_n-\alpha A_n=\nabla \theta_n^0-\alpha A_n^0,
\end{equation}
 and in particular \eqref{e:thetacondition} is zero outside of $\bigcup_i\Omega_i$.
 We now prove \eqref{e:gammalimsup}. Notice that, because of \eqref{e:thetacondition} and the definitions \eqref{e:rho},\eqref{e:b}, \eqref{e:aglobal} and \eqref{e:theta}, we have \begin{equation*}
    \mathcal{E}_{\eps_n}(u_n,A_n,T^c)=0.
\end{equation*} We just have to focus on the squares $Q_{\frac{\eps_n}{\eps_0}}(x_{i,j})$ composing the sets $\Omega_i$ and on the ``corners" (that is $T\backslash\bigcup_i \Omega_i$). By a change of variable and using Proposition \ref{prop:buildingblock},  we have for every $i=1\dots N$ and $j=1\dots N_i$
\begin{equation}
\label{e:squares}
    \mathcal{E}_{\eps_n}(u_n,A_n,Q_{\frac{\eps_n}{\eps_0}}(x_{i,j}))=\frac{\eps_n}{\eps_0}\mathcal{E}_{\eps_0}(u_0,A_0,Q_1)\leq \frac{\eps
    _n}{\eps_0}(\sigma_\kappa+\xi).
\end{equation}
We now show that the energy near the corners tends to $0$ as $\eps_n\rightarrow 0$. By \eqref{e:1denergycontrol} and since \eqref{e:thetacondition} and $B_n$ are zero on $T\backslash\bigcup_i \Omega_i$, we have
\begin{equation}\label{e:corners}
    \mathcal{E}_{\eps_n}(u_n,A_n, T\backslash\bigcup_i \Omega_i)=\int_{ T\backslash\bigcup_i \Omega_i}\eps_n\abs{\nabla\rho_n}^2+\frac{1}{2\eps_n}(1-\rho_n^{2})^{2}\,dx\lesssim \frac{1}{\eps_n}\abs{ T\backslash\bigcup_i \Omega_i}\stackrel{\eqref{e:smallangles}}{\lesssim} \eps_n.
\end{equation}
Putting \eqref{e:squares} and \eqref{e:corners} together we get,
\begin{equation*}
    \mathcal{E}_{\eps_n}(u_n,A_n)=\sum_{i,j}\mathcal{E}_{\eps_n}(u_n,A_n,Q_{\frac{\eps_n}{\eps_0}}(x_{i,j}))+\mathcal{E}_{\eps_n}(u_n,A_n, T\backslash\bigcup_i \Omega_i)\leq \sum_iN_i\frac{\eps
    _n}{\eps_0}(\sigma_\kappa+\xi)+C\eps_n.
\end{equation*}
Since $N_i\eps_n/\eps_0\rightarrow L_i$ as $\eps_n\rightarrow 0$, taking the $\limsup$ in the previous inequality yields 
\begin{equation*}
    \limsup_{\eps_n} \mathcal{E}_{\eps_n}(u_n,A_n)\leq  \sum_iL_i(\sigma_\kappa+\xi)=P(\{\rho=0\},Q_1)(\sigma_\kappa+\xi).
\end{equation*}
The conclusion then follows by a diagonal procedure with respect to $\xi$.

\quad \emph{\textbf{Case 2: $E$ connected closed polyhedral set. } }
We now assume that $E$ is a connected closed polyhedral set. We can repeat the constructions of Case 1 for $T$, $\rho_{n},B_{n},A_n,A^0_{n}$, $\theta^0_{n}$ and $A_n$, as in \eqref{e:epsneigh},  \eqref{e:rho}, \eqref{e:b},\eqref{e:aglobal}, \eqref{e:a_0} and \eqref{e:thetazero} respectively, since they do not use the fact that $E$ is simply connected. We can also choose the value $\zeta$ as in \eqref{e:hest} so that \eqref{e:massconstr} holds.\\
We are left with the definition of an appropriate phase $\theta_n$. Since $E$ is not simply connected, $E^c\cup T$ is not connected in general. Let us call $S>0$ the number of connected components of $E^c\cup T$ and $F_k, \, k=1\dots S$ each component. Since $E$ is connected, at most one connected component of $E^c\cup T$ is not simply connected and encloses $E$. Without loss of generality, we assume that $F_1$ encloses $E$ and $F_k,\, k=2\dots S$ are simply connected. We now define $\theta_n$ separately on each connected component of  $E^c\cup T$.  For every $k=1\dots S$, we define $\theta_n(x)$ for $x\in F_k\cap \{\rho>0\}$ as 
\begin{equation*}
    \theta_n(x)=\theta^0_n(x)+\alpha\int_{\Gamma_{x,x_{0}^k}}(A_n-A^0_n)\cdot\,dl
\end{equation*}
where $x_{0}^k\in F_k$ is arbitrary and $\Gamma_{x,x_{0}^k}$ is any curve contained in $F_k$ connecting $x$ to $x_{0}^k$. For $k=1$, the same computations as \eqref{e:quantization} give that $\theta_n$ is defined up to a multiple of $2\pi$. Then, if $k>1$, the set $F_k$ is simply connected. Since $A_n-A_n^0$ is irrotational, the definition of $\theta_n$ does not depend on the choice of the curve $\Gamma_{x,x_{0}^k}$. This means that as before $u_n=\rho_ne^{i\theta_n}$ is uniquely defined, and the same computations as Case 1 apply for the proof of the $\Gamma$-limsup inequality.

\quad \emph{\textbf{Case 3: $E$ closed polyhedral set.} }
Let us finally assume that $E\subset \tor$ is any arbitrary closed polyhedral set. In particular $E$ has a finite number of sides and thus a finite number $M$ of connected components. We denote $E_j$ with $j=1\dots M$ these sets. We now consider for $j=1 \dots M$ values $\zeta_j \in (-\eps_n/(2\eps_0),\eps_n/(2\eps_0))$ to be chosen later. We set $T_j=T(\eps_n,\zeta_j, E_j)$ and $T=~\cup_{j=1}^M T_j$. For $\eps_n$ small enough we can suppose that the sets $T^j$ are disjoint. For each set $T_j$ and value $\zeta_j$ consider the constructions of Case 1 \eqref{e:rho}, \eqref{e:b}, \eqref{e:a_0} and \eqref{e:thetazero} which we denote by $\rho_{n,j},B_{n,j},A^0_{n,j}$ and $\theta^0_{n,j}$. We set
\begin{equation*}
\rho_n(x)=
    \begin{cases}
        0 &\text{in}  \ E_j\backslash  T_j,\\
        \rho_j(x) &\text{in} \ T_j,\\
        1 &\text{elsewhere}
    \end{cases}  
\end{equation*}
and
\begin{equation*}
B_n(x)=
    \begin{cases}
        \frac{1}{\sqrt{2}} &\text{in}  \ E_j\backslash  T_j,\\
        B_j(x) &\text{in} \ T_j,\\
        0 &\text{elsewhere}.
    \end{cases}  
\end{equation*}
 Arguing as in \eqref{e:hest}, we have for $j=1\dots M$
\begin{equation*}
    \int_{E_j\cup T^j}B_n\,dx=\frac{1}{\sqrt{2}}\abs{E_j}-P(E_j)\left[\frac{1}{\sqrt{2}}\left(\zeta_j+\frac{\eps_n}{2\eps_0}\right)-\frac{\eps_n}{\eps_0}\int_{Q_1}B_0\,dx\right]+o(\eps_n).
\end{equation*}
We now choose in a similar way as before the values $\zeta_j\in (-\eps_n/(2\eps_0),\eps_n/(2\eps_0))$ such that for all $j=1\dots M-1$ we have (recall that $\alpha=\kappa^-1\eps_n^{-2}$)  
\begin{equation*}
\int_{E^j\cup T^j_n}B_n\,dx =\frac{2\pi}{\alpha}\left\lfloor\frac{\alpha}{2\pi}\frac{\abs{E_j}}{\sqrt{2}}\right\rfloor,
\end{equation*} so that 
\begin{equation}
\label{e:localquantization}
     \alpha\int_{E^j\cup T^j_n}B_n\,dx\in 2 \pi \Z.
\end{equation}
Moreover, we choose $\zeta_M\in (-\eps_n/(2\eps_0),\eps_n/(2\eps_0))$ such that
\begin{equation*}
\int_{\tor} B_n\,dx=b_{\text{ext}}/\kappa.    
\end{equation*}
The assumption $\eps_n^{-2}\kappa^{-2}b_{\text{ext}}\in 2\pi\Z$ then ensures that \eqref{e:localquantization} holds also for $j=M$. As before, we define the auxiliary field for  $x\in E^c\cup T$
\begin{equation*}
A_n^0(x)=
    \begin{cases}
        A_{n,j}^0(x) &\text{in} \ T_j,\\
        0 &\text{elsewhere}
    \end{cases}  
\end{equation*}
and for $x\in \{\rho_n>0\}$
\begin{equation*}
\theta_n^0(x)=
    \begin{cases}
        \theta_{n,j}^0(x) &\text{in} \ T_j,\\
        0 &\text{elsewhere}.
    \end{cases}  
\end{equation*}  
We then define $A_n$ as a global vector potential on $Q_1$ as in \eqref{e:aglobal}. Now consider the set $E^c\cup T$. Let us call $S>0$ the number of its connected components and denote by $F_k$, $k=1\dots S$ each component. As in Case 2, we define the phase $\theta_n$ locally on each connected component of $E^c\cup T$. For every $k=1\dots S$, we define $\theta_n(x)$ for $x\in F_k\cap \{\rho>0\}$ as 
\begin{equation*}
    \theta_n(x)=\theta^0_n(x)+\alpha\int_{\Gamma_{x,x_{0}^k}}(A_n-A^0_n)\cdot\,dl
\end{equation*}
where $x_0^k\in F_k$ is arbitrary and $\Gamma_{x,x_0^k}$ is any curve contained in $F_k$ connecting $x$ to $x_0^k$. Differently to Case 2, each $F_k$ may now enclose multiple connected components of $E$, see Figure \ref{fig:components}.
\begin{figure}[h]
    \centering
     \includegraphics[scale = 1.2]{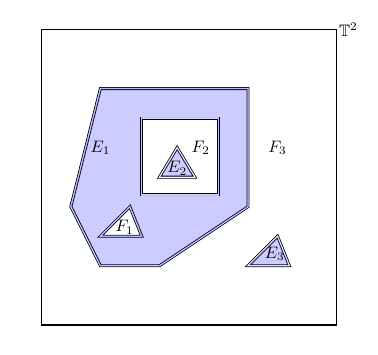} 
    \caption{Example of a polyhedral set $E$ (in blue) such that both $E$ and $E^c$ have multiple connected components, $E_1$, $E_2$, $E_3$, and $F_1$, $F_2$, $F_3$ respectively.  The set $T$ is represented by the double line around $\partial E$. The phase $\theta_n$ is defined locally on each $F_k,\,  k=1\dots 3$. On any $F_k, \, k=1\dots 3$, to ensure that the phase is defined up to a multiple of $2\pi$, each connected component of $E$ which is enclosed by $F_k$ is subject to the flux quantization condition \eqref{e:localquantization}. In particular, in this example $F_3$ forces a quantization condition on $E_1, E_2$ and $E_3$, while $F_2$ forces a quantization condition on $E_2$. Since $F_1$ is simply connected, its phase is already quantized and no quantization condition is enforced. }
    \label{fig:components}
\end{figure} 
To prove that $\theta_n$ is defined up to a factor of $2\pi$, for any $k=1\dots S$ consider a closed curve $\Gamma$ in $F_k$. With the same computation as \eqref{e:quantization} we get
    \begin{equation}
    \label{e:multiplepoles}
    \alpha\int_\Gamma (A-A_0)\cdot\,dl = \alpha\sum_ja_j\int_{E^j\cup T^j_n}B_n\,dx,
\end{equation}
for some $a_j\in\{-1,1 \}$, where the sum on the right hand side of \eqref{e:multiplepoles} is performed on all $j$ such that $E_j$ is enclosed by $\Gamma$.  By \eqref{e:localquantization} all the terms in the sum are multiples of $2\pi$, hence also the sum is a multiple of $2\pi$. Thus,  $u_n=\rho_ne^{i\theta_n}$ is uniquely defined. To conclude, the same computations as in Case 1 apply for the proof of the $\Gamma$-limsup inequality.

\end{proof}
\section*{Acknowledgement}
 A.C. is partially supported by the European Union's Horizon 2020 research and innovation program under the Marie Sklodowska-Curie grant agreement No 945332. Part of this research was supported by the ANR project STOIQUES. This work was supported by a public grant from the Fondation Mathématique Jacques Hadamard.
\begin{flushright}
    \includegraphics[height=.9cm]{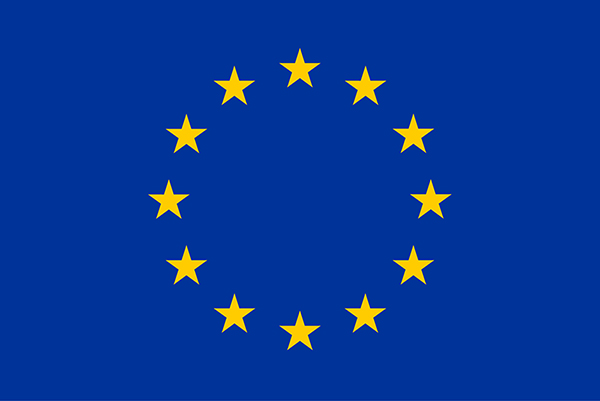}
\end{flushright}
\appendix
\section{Derivation of the functional from the Ginzburg-Landau Model}
\label{sec:derivation}
Let us derive the functional \eqref{e:GL2D} from the full Ginzburg-Landau functional. We discuss the derivation only at an heuristic level. Given a sample $\Omega\subset \R^3$, the Ginzburg-Landau functional is 
\begin{equation*}
 GL(u,A)=\int_{\Omega} |\nabla_A u|^2 +\frac{\kappa^2}{2} (1-|u|^2)^2 dx +\int_{\R^3} |\nabla\times A- \Bext|^2 dx,
\end{equation*}
where $u:\Omega\rightarrow \C$ is the order parameter representing the state of the material, $A:\R^{3}\rightarrow\R^{3}$ is the magnetic vector potential of the magnetic field $B=\nabla \times A$, $\nabla_{A}u=\nabla u-iAu$ denotes the covariant derivative of $u$ and $\Bext$ is the external magnetic field. The parameter $\kappa\in\R$ is the Ginzburg-Landau parameter. Since we focus on type-I superconductors, we assume $0<\kappa<1/\sqrt{2}$. 

We consider a sample $\Omega = \tor_L\times (-T,T) = (\R/L\Z)^2 \times (-T,T)$ of sides $L>0$ and thickness $T>0$ and assume lateral periodic boundary conditions. Moreover, we assume that the external magnetic field is uniform and perpendicular to the sample, that is $\Bext=(0,0,\bext)$.

We now write
\begin{equation*}
    \mathcal{D}_{A}^{3}u=(\nabla_A u)_{2}-i(\nabla_{A}u)_{1}
\end{equation*}
and, for any vector field $v:\R^3\rightarrow\R^3$, we write $v'=(v_1,v_2,0)$. Arguing as in \cite[Lemma 2.3]{COS}, we have
\begin{equation*}
    GL(u,A)=2LT^2(\kappa \sqrt{2}\bext-\bext^2)+F(u,A),
\end{equation*}
where the first term, the bulk energy of the Ginzburg-Landau functional, does not depend on $u$ and $A$, and
\begin{multline}
\label{e:nonbulkGL}
    F(u,A)=\int_{\tor_L\times (-T,T)}(1-\kappa\sqrt{2})\abs{\nabla_{A}'u}^{2}+\kappa\sqrt{2}\abs{\AderD{u}{A}}^2+\abs{(\nabla_A u)_3}^2\,dx\\
    +\int_{\tor_L\times (-T,T)}\left(B_3-\frac{\kappa}{\sqrt{2}}(1-\abs{u}^{2})\right)^2\,dx \\+
     \int_{\tor_L\times \R}\abs{B'}^2\,dx+\int_{\tor_L\times (\R/(-T,T))}(B_3-\bext)^2\,dx.
\end{multline}
Since $\diverg B=0$, the periodicity in the first two directions implies that for any $x_3\in(-T,T)$
\begin{equation*}
    \int_{\tor_L\times \{x_3\}} B_3\,dx=L^2\bext,
\end{equation*}
see also \cite[Lemma 2.2]{COS}.
Let us now take $T=+ \infty$. By symmetry, this implies that $u$ and $A$ should not depend on $x_3$ and $B$ should be parallel to $e_3$. We can thus impose $u(x)=u(x')$ and $A(x)=(A_1(x'),A_2(x'),0)$, which implies $B(x)=(0,0,B_3(x'))$. We can now consider the energy density on a slice $\tor_L\times \{x_3\}$, which with abuse of notation we still call $F(u,A)$. Plugging in \eqref{e:nonbulkGL} the assumptions on $u$ and $A$ we get
\begin{equation*}
    F(u,A)=\int_{\tor_L\times\{x_3\}}(1-\kappa\sqrt{2})\abs{\nabla_{A}'u}^{2}+\kappa\sqrt{2}\abs{\AderD{u}{A}}^2+ \left(B_3-\frac{\kappa}{\sqrt{2}}(1-\abs{u}^{2})\right)^2\,dx.
\end{equation*}
Notice that the last two terms in \eqref{e:nonbulkGL} disappear by our assumptions on $B$ and $T$. We can see that $F(u,A)$ resembles a Modica-Mortola type functional. Since we are interested in the limit behavior of the functional, and since we expect the functional to behave like a perimeter term in 2D, we rescale it by the typical length $L$. Dropping the prime superscripts and calling $B=B_3$ we get
\begin{equation}
    \label{e:GL2Dor}
    \mathcal{E}(u,A)=\frac{1}{L}F(u,A)=\frac{1}{L}\int_{\tor_L}(1-\kappa\sqrt{2})\abs{\Ader{u}{A}}^{2}+\kappa\sqrt{2}\abs{\AderD{u}{A}}^2+
    \left(B-\frac{\kappa}{\sqrt{2}}(1-\abs{u}^{2})\right)^2\,dx
\end{equation}
under the constraint 
\begin{equation*}
    \frac{1}{L^2}\int_{\tor_L}B\,dx=b_{\text{ext}}.
\end{equation*}
Finally, we consider the rescaling $x=L\hat{x}$, $A(x)=\kappa L \hat{A}(\Hat{x})$, $u(x)=\hat{u}(\hat{x})$  and we set $\eps=(\kappa L)^{-1}$, so that $\Ader{u}{A}=\left(L\right)^{-1}\hat{\nabla}_{\kappa^{-1}\eps^{-2}\hat{A}}\hat{u}$ and $B=\kappa\hat{B}$. We find (dropping hats)
\begin{equation*}
    \mathcal{E}_\eps(u,A)=\int_{\tor}\eps((1-\kappa\sqrt{2})\abs{\Ader{u}{\kappa^{-1}\eps^{-2}A}}^{2}+\kappa\sqrt{2}\abs{\AderD{u}{\kappa^{-1}\eps^{-2} A}}^2)+\frac{1}{\eps}\left(B-\frac{1}{\sqrt{2}}(1-\abs{u}^{2})\right)^{2}\,dx,
\end{equation*}
with \begin{equation*}
    \int_{\tor}B\,dx=\frac{b_{\text{ext}}}{\kappa}.
\end{equation*}

\printbibliography

\medskip
\small
\begin{flushright}
\noindent \verb"acosenza@math.univ-paris-diderot.fr"\\
Universit\'e Paris Cit\'e and Sorbonne Universit\'e,\\
CNRS, Laboratoire Jacques-Louis Lions (LJLL),\\
F-75006 Paris, France
\end{flushright}
\medskip
\small
\begin{flushright}
\noindent \verb"michael.goldman@cnrs.fr"\\
CMAP, CNRS, \'Ecole polytechnique, Institut Polytechnique de Paris,\\ 
91120 Palaiseau, France
\end{flushright}
\medskip
\begin{flushright}
\noindent \verb"azilio@math.univ-paris-diderot.fr"\\
Universit\'e Paris Cit\'e and Sorbonne Universit\'e,\\
CNRS, Laboratoire Jacques-Louis Lions (LJLL),\\
F-75006 Paris, France
\end{flushright}

\end{document}